\documentclass[]{interact}

\usepackage{epstopdf}% To incorporate .eps illustrations using PDFLaTeX, etc.
\usepackage{xcolor}
\usepackage{subcaption}
\usepackage{algorithm}
\usepackage[noend]{algpseudocode}
\usepackage[font=scriptsize]{caption}
\usepackage[labelformat=simple]{caption}

\usepackage[natbibapa,nodoi]{apacite}
\theoremstyle{plain}% Theorem-like structures provided by amsthm.sty

\theoremstyle{definition}

\theoremstyle{remark}

\begin{document}

%\articletype{ARTICLE TEMPLATE}% Specify the article type or omit as appropriate

\title{State Discretization for Continuous-State MDPs in Infectious Disease Control}
%\maketitle
\author{
\name{Suyanpeng Zhang $^a$\thanks{CONTACT Suyanpeng Zhang. Email: suyanpen@usc.edu} and Sze-chuan Suen $^a$}
\affil{\textsuperscript{a}Daniel J. Epstein Department of Industrial and Systems Engineering, Viterbi School of Engineering, University of Southern California, Los Angeles, CA, USA}}

\maketitle

\begin{abstract}
Repeated decision-making problems may arise in the health policy context, such as infectious disease control for COVID-19 and other epidemics. These problems may sometimes be effectively solved using Markov decision processes (MDPs). However, the continuous or large state space of such problems for capturing infectious disease prevalence renders it difficult to implement tractable MDPs to identify the optimal disease control policy over time. We therefore develop an algorithm for discretizing continuous states for approximate MDP solutions in this context. We benchmark performance against a uniform discretization using both a synthetic example and an example of COVID-19 in Los Angeles County.
\end{abstract}


\section{Introduction} \label{intro}
Public health officials often need to determine the optimal population health intervention policy over time even as the state and trajectory of disease are driven by complex dynamics. Many of these problems require making policy decisions sequentially over time, where the state may be represented using a continuous measure (e.g., the proportion of the population that is infected). For instance, during the COVID-19 pandemic, decision-makers needed to repeatedly set the start and end times of lockdowns that limited travel and interactions between individuals to reduce transmission without fully understanding the exact transmissibility of COVID-19. This sequential decision-making problem appears repeatedly in infectious disease control problems, as evidenced by prior literature on similar problems \citep{kaplan1996general,blower2002health,talbot2005influenza,zhang2011simulation,matrajt2021optimizing,fu2022optimal}. Such problems often take into account underlying disease dynamics, which are uncertain or depend on a variety of complex social and biological factors.

A difficulty in solving repeated decision making problems for infectious disease control is the complexity of infectious disease dynamics, which are typically represented using compartmental models and simulation-based models \citep{brauer2008compartmental,kopec2010validation}. Such models are difficult to use for repeated decision-making problems as one often needs to evaluate the model repeatedly to identify an optimal policy for disease control, which may require a significant investment of computational time, as there is no closed-form solution. 

While there are sophisticated means to identify optimal policies, these techniques have their own challenges. For instance, the maximum principal approach
\citep{goenka2014infectious,pontryagin2018mathematical,piguillem2022optimal} offers a solution framework for optimal control issues under differential equation systems. However, its application becomes increasingly challenging with a large number of states or policies. Such expansion complicates both the Hamiltonian and the differential equations system, thereby rendering the process of deriving analytical or numerical solutions complicated and time-consuming. Moreover, it is difficult to find the optimal solution when the problem is non-convex. Simulation optimization, which can handle complex systems, has also been used in disease contexts \citep{carson1997simulation}. However, this can also be computationally expensive and time-consuming. Furthermore, the quality of the solution highly depends on the search space and the heuristic function chosen, presenting challenges to its practical application. This problem can be formulated as a dynamic programming problem \citep{calvia2023simple}, but the continuous or large state space can create difficulties. Furthermore, while the infectious disease control problem can be formulated as a mixed-integer programming problem, the inherently non-linear nature of most disease dynamics—such as the compartmental model—introduces non-linear constraints into the formulation. Consequently, the problem becomes extremely challenging to solve \citep{bertsimas2022locate}, which limits its generalizability, especially as the disease dynamics become more complex.

Markov decision processes (MDPs) are also a commonly used method for repeated decision making problems. MDPs allow for state transitions, which can be used to describe changes in disease/health states over time and allow for repeated decisions over time. Given current computing innovations, many MDPs of useful size can be solved effectively using algorithms such as backward induction, value iteration, policy iteration, etc. MDPs can also be efficiently solved with non-convex problems. 

However, incorporating dynamics from compartmental models and simulations into an MDP framework is challenging because disease models often use a continuous or large number of possible states (as the state usually represents a proportion of the whole population in certain statuses like infected, recovered, and hospitalized). Having a continuous state space makes the MDP problem difficult to solve since traditional MDP solution methods may then not work even for a short time horizon due to state-space explosion issues. For example, backward induction need $|S|^2|A||T-1|$ multiplications. In the case of value iteration, each iteration carries a complexity of $\mathcal{O}(|S|^2|A|)$. In the case of policy iteration, each iteration carries a complexity of $\mathcal{O}(|S|^3+|S|^2|A|)$, and modified policy iteration requires $\mathcal{O}(k|S|^2+|S|^2|A|)$ per iteration \citep{Puterman}. For this reason, many traditional MDP studies in the healthcare field focus on finite-state decision-making problems like monitoring, treatment initiation, and disease testing and diagnosis \citep{Claude,David1985,Ahn,Hu,	MAGNI2000237,Alagoz2004, Alagoz2007,Jennifer,Maillart2008,Shechter2008,denton2009,Chhatwal2010,Kurt2010,Alagoz2013,MASON2014727,Capan2017,Liu2017,zhang2021early, Suen2018}. Therefore, finding a good state discretization method that translates infectious disease dynamics onto a limited number of states improve computational efficiency and potentially widen the scope of MDP applications, particularly in the infectious disease space. 

Uniform discretization is a traditional way of addressing continuous state problems. However, this methodology is suboptimal for addressing infectious disease control challenges. The heterogeneity in state visit frequencies—wherein some states (with extremely high prevalence) may remain unvisited and others (with lower prevalence) might be visited more frequently—renders uniform discretization inefficient. This approach may result in the overuse of discretization regions towards states that are less likely to be visited and an inadequate number of discretization regions for those with higher probabilities of being reached. How can we find a better way of discretizing the state space to closely represent the changes in health systems/disease? While many works have used various discretization methods to reduce state spaces \citep{lovejoy1991computationally,sandikcci2013alleviating}, we take a novel approach that treats the state discretization problem as an optimization problem. This allows us to find the discretization that will provide a smaller discretized region in more likely visited states for a more accurate description of the true dynamics.

We will explore the above state discretization in the context of a disease control problem where states are used to describe the disease dynamics over a population, actions are implemented to prevent disease spread (lockdown, social distancing, face masks, and so on). States are assumed to be fully observable at each time period. Under this framework, we find a better way of discretizing states such that the discretized state space serves as a good proxy of the original state space.
This paper addresses the challenge of formulating infectious disease control problems as MDPs by proposing a new algorithm for non-uniform state discretization that enables the discrete representation of infinite state spaces.

\subsection{Contributions}
We make several contributions in this study. We provide a novel algorithm for defining a non-uniform, discrete state space for infectious disease control problems that well approximates the original continuous state dynamics. Our algorithm exploits the likelihood of each state being visited in the system to more efficiently capture the transitions between states. Defining a discrete set of states from an originally continuous system allows us to incorporate infectious disease dynamics within frameworks that are better suited for discrete state spaces, such as MDPs. Finally, we demonstrate that our state space discretization allows for more accurate MDP outcomes through two numerical examples, one using a classic SIR compartmental model and one using COVID-19 model of Los Angeles County.

The remainder of this paper is organized as follows: we review the related literature in Section~\ref{liter_discretization}, present the problem setup in Section~\ref{problem_setup_discretization}, and provide the algorithms in Section~\ref{algo}. The numerical example is shown in Section~\ref{numeric_discretization}. In Section~\ref{conclusions_discretization}, we conclude.

\section{Literature Review}\label{liter_discretization}
\subsection{Markov Decision Processes in Healthcare Applications}

MDPs have a rich history in the field of operations research, with wide range of applications such as inventory management \citep{giannoccaro2002inventory}, portfolio management \citep{bauerle2009mdp}, production and storage optimization \citep{arruda2008stability}, and various others. Extensive research has been conducted to solve and understand the structure of MDPs, with notable contributions from works such as \citet{Puterman} and \citet{topkis2011supermodularity}. 

MDPs also find widespread application in the field of healthcare. They offer valuable insights and solutions to various health-related issues, including scheduling \citep{agrawal2023preference}, screening \citep{Maillart2008,Chhatwal2010,Alagoz2013,mcnealey2023optimizing,shen2023advance}, sequential disease testing \citep{arruda2019optimal,singh2020partially}, treatment initiation \citep{Shechter2008,Liu2017,otero2023monitoring}, and organ transplantation \citep{sandikcci2008estimating,sandikcci2013alleviating,zhang2021early,zhang2024quantifying}. For instance, patients in different age groups with risks of breast cancer may need personalized mammography exam frequencies \citep{Alagoz2013}, or, in another example, a patient with organ failure may be presented at different states with organ transplant options that vary in their compatibility with the patient. The patient may face the decision to either wait for a better match or accept an offered organ as their own survival probability decreases over time \citep{zhang2021early,zhang2024quantifying}. For a more extensive exploration of MDPs in healthcare, refer to the comprehensive reviews by \citet{Schaefer2004}, \citet{Alagoz2010MarkovDP}, and \citet{Sonnenberg}. Although MDPs are widely used in healthcare applications, most of these consider finite-state decision-making. Constructing an MDP for infectious disease control problems with repeated decisions is challenging, especially when the state space for such problems is continuous.

\subsection{Solving Continuous State MDP}\label{liter:discretization}
As previously discussed, an infinite or continuous state space is a major challenge when formulating MDPs. Several methods have been proposed to address this problem. In \citet{li2005lazy}, a discretization-free approach (modified value iteration with lazy approximation) is introduced, but it requires highly precise piece-wise constant approximations of both the value and transition functions. In \citet{munos2002variable}, different criteria for discretizing state and time space non-uniformly are discussed. These methods involve evaluating values or policies using dynamic programming; however, some of these methods raise computational concerns for problems with continuous or large numbers of states. \citet{zhou2010solving} used Monte-Carlo simulation to approximate the belief state using a finite number of particles on a discretized grid mesh. However, the study does not provide guidance on the construction of the grid mesh. \citet{brooks2006parametric} proposed a parametric method to uniformly discretize a continuous state space over a lower dimensional parametric space. However, since prior knowledge of the distribution is required, MDPs for infectious disease control would be difficult to solve in this manner. One remedy is to solve the MDP formulation by truncation and discretization of the state \citep{boucherie2017markov}. Researchers have used various methods to achieve this. For example, \citet{sandikcci2013alleviating} used fixed-resolution, non-uniform grids to discretize the belief state and approximate the optimal policy for a partially observable MDP (POMDP) model. \citet{lovejoy1991computationally} used fixed or uniform grids to approximate the solution of the POMDP. However, using uniform or pre-defined discretization regions (which requires domain knowledge) may not always be appropriate, particularly for infectious disease control problems where disease spread is subject to substantial changes across different policy scenarios. In such cases, a more effective discretization algorithm is needed to enable the computation of the optimal policy.

\subsection{Modeling Disease Dynamics}\label{population_model}
To identify the optimal policy for an infectious disease control problem, it is necessary to have a model for describing the disease dynamics. For instance, during an emerging pandemic, how would disease transmission change if the government imposed a 1-month lockdown? How would it change if the government imposed a 3-month lockdown instead? Different policies may change the patterns of disease transmission and thus change the proportion of infections in total. To efficiently avert infections, these different possibilities need to be evaluated to understand the resultant health and cost outcomes. Multiple methods are available for assessing the impact of different policies on a specific population. 

One common method to model disease dynamics is to use compartmental models based on differential equations \citep{brauer2008compartmental,kermack1991contributions,kermack1991contributions2,kermack1991contributions3}. A compartmental model uses a mathematical framework to provide insights into the mechanism that affect the transmission and progression of disease. This framework partitions the population into different health or treatment states (compartments). For instance, each compartment represents a specific stage of the infectious disease (e.g., susceptible, infected, recovered), and proportions of the population move between compartments described by differential equations at certain rates. This model is fundamental in epidemiology for understanding the spread of diseases and evaluating the potential impact of public health interventions. For example, compartmental models can compare the effectiveness of wearing masks and social distancing during the COVID-19 pandemic \citep{grimm2021extensions,kai2020universal}. \citet{long2018spatial} use a classical compartmental model to assist with the decision of allocating resources during the 2014 Ebola outbreak in Africa. In Section~\ref{numeric_discretization}, we consider a classic Susceptible-Infected-Recovered (SIR) epidemic model, which has been extensively used in the epidemiological literature \citep{beckley2013modeling,harko2014exact,kroger2020analytical}.

Another method of evaluating disease dynamics is to use simulation models, which can be used to track transmission, progression, and behavior as well as policy outcomes. For instance, simulation models can be employed to examine the cost-effectiveness of screening recommendations for positive-HIV men who have sex with men (MSM) \citep{tuite2014cost}, as well as to study the effectiveness of different disease control strategies for tuberculosis (TB) in India \citep{suen2014disease}. Although these methods indeed capture the dynamics of complicated diseases, they are unable to compute dynamic policies effectively as $m^t$ evaluations are usually needed when there are $m$ possible interventions and $t$ decision epochs. Therefore, it is beneficial to find alternative effective ways of identifying the optimal policy for infectious disease control. In our paper, we consider a discrete-state MDP framework that takes advantage of its effective solution methods with underlying disease dynamics estimated from traditional disease models such as compartmental and simulation models.

To model this problem as a discrete-state MDP, we also need to define a transition function to describe the probability of transitioning between the states. Several existing techniques can be used to construct this function. For instance, \citet{yaesoubi2011generalized} proposed a way to compute transition probabilities given a system of ODEs. In another example, \citet{mishalani2002computation} proposed a method of developing transition probabilities from a stochastic duration model based on the hazard rate function. However, these methods are computationally intensive, which limits their usage to problems with small populations or disease models with special structures. 

\section{Problem Setup}\label{problem_setup_discretization}

\begin{table}[ht]
    \centering
    \begin{tabular}{|c|c|}
    \hline
    $T$ & The set of all decision epochs \\ \hline
    $A$ & The set of all policy interventions\\ \hline
    $\mathcal{X}$ & The set of all state representations\\ \hline
    $\bar{\mathcal{X}}$ & The set of all discretized state representations\\ \hline
    $\{\pi_t\}$ & Policy intervention at decision epoch $t$ \\ \hline
    $G_d$ & Discretization vector for component $d$ \\ \hline
        $\{\mathbf{X}_t\}$ & Observed disease trajectory \\ \hline
        $\{\bar{\mathbf{X}}_t\}$ & Discretized disease trajectory \\ \hline
        $\{\tilde{\mathbf{X}}_t\}$ & Markovian disease trajectory computed using transition probabilities\\ \hline 
        $f(\mathbf{X}_t,\pi_t)$ & Transition function\\ \hline
        $\bar{f}(\mathbf{X}_t,\pi_t,G)$ & Discretized transition function\\ \hline
        $P(\pi_t)$ & Transition probability for policy $\pi_t$\\ \hline
        $\lambda$ & Discount factor\\ \hline
        $r(\mathbf{X}_t,\pi_t)$ & Reward function\\ \hline
    \end{tabular}
    \caption{Table of notation.}
    \label{tab:notations}
\end{table}
The notation used in this paper is as follows. We denote $\mathbf{X}_t\in \mathcal{X}$ as the state of the epidemic at time $t$. $\mathbf{X}_t=[X_{1t},X_{2t},...,X_{nt}]$ has $n$ components where each represents the proportion of the population in the compartment (e.g., for a SIR model, $n=3$). For example, $\mathbf{X}_t = [X_{St},X_{It},X_{Rt}]\in [0,1]^3$ can describe the proportion of the population in susceptible (S), infected (I), and recovered (R) compartments at time $t$ for a SIR model. We denote $\mathbf{X}_0$ as the initial state and we assume it follows an initial distribution $\Omega$. We use $\{\mathbf{X}_t\}=(\mathbf{X}_0,...,\mathbf{X}_N)$ to denote the disease trajectory.

In this paper, we focus on the finite horizon problem. Let $T = \{1, ..., N \}$ be the set of possible decision epochs for the problem. $A= \{1, ..., |A| \}$ is the set of possible policy interventions for the problem. We assume a small, finite number of actions/policies (e.g., lockdown versus no lockdown). We denote $\pi_t\in A$ as the policy intervention at time $t$. 

We consider a model denoted by $f(\mathbf{X}_t,\pi_t)=\mathbf{X}_{t+1}$ that describes the disease dynamics across time epochs $t$. This function $f(\mathbf{X}_t,\pi_t)$ can consider disease progression, transmission over time, mortality, and interventions. Generally, $f(\mathbf{X}_t,\pi_t)$ takes the state of the system and policy intervention as an input and then returns the state in the next period. We assume that $f(\mathbf{X}_t,\pi_t)$ is time-homogeneous for simplicity (if time-inhomogeneous dynamics are desired, our methods can be easily extended). 

The cost in state $\mathbf{X}_t\in \mathcal{X}$ and taking action $\pi_t\in A$ for $t\in T$ in the infectious disease control problem is denoted using $r(\mathbf{X}_t,\pi_t)$. This cost can be dependent on health outcomes (e.g., number of infected, total vaccinated population, etc.) as well as other factors (financial cost, economic burden, etc.). We let $\lambda$ denote the discount factor.

Given the transition function $f(\mathbf{X}_t,\pi_t)$ and the cost function $r(\mathbf{X}_t,\pi_t)$, we have the following optimization formulation for our repeated decision-making disease control problem:
\begin{align}
    \min_{\pi_0,...,\pi_{N-1}} & \ \ \sum_{t=0}^N \lambda^t r(\mathbf{X}_t,\pi_t)|\mathbf{X}_0\\
    s.t. & \ \ \mathbf{X}_t = f(\mathbf{X}_{t-1},\pi_{t-1})
\end{align}

In the above problem, the objective is to find a sequence of actions $\{\pi_0,...,\pi_{N-1}\}$ that minimizes the total discounted cost function $r(\mathbf{X}_t,\pi_t)$ over states $\mathbf{X}_t$ for the whole $N$-period time horizon given a known initial state $\mathbf{X}_0$. For example, $\mathbf{X}_t$ can represent the proportion of individuals in each COVID-19-related health stage at time $t$, and let $r(\mathbf{X}_t,\pi_t)$ compute the proportion of people dead from COVID-19 at time $t$. If $\pi_{t}$ denotes the policy intervention (lockdown or not) at time $t$, then $f(\mathbf{X}_{t-1},\pi_{t-1})$ could be a system of difference equations that describes the population flow across different health stages. Our objective in this problem then is to find the optimal policy intervention at each time $t$ that minimizes the total cost within $N$ periods. 

There are challenges to solving the above formulation using traditional MDP solution methods (e.g., backward induction, value iteration, policy iteration, etc.) as this formulation usually contains constraints with non-linear dynamics on a continuous state space. These solution methods require a finite number of states for effective evaluation. Moreover, the function $f(\mathbf{X}_t,\pi_t)$ may not be expressed as transition probabilities from state to state, while many traditional MDP solution methods use transition probability matrices to allow for the modeling of uncertainty and variability in decision-making processes. 

\begin{figure}
    \centering
    \includegraphics[width=4in]{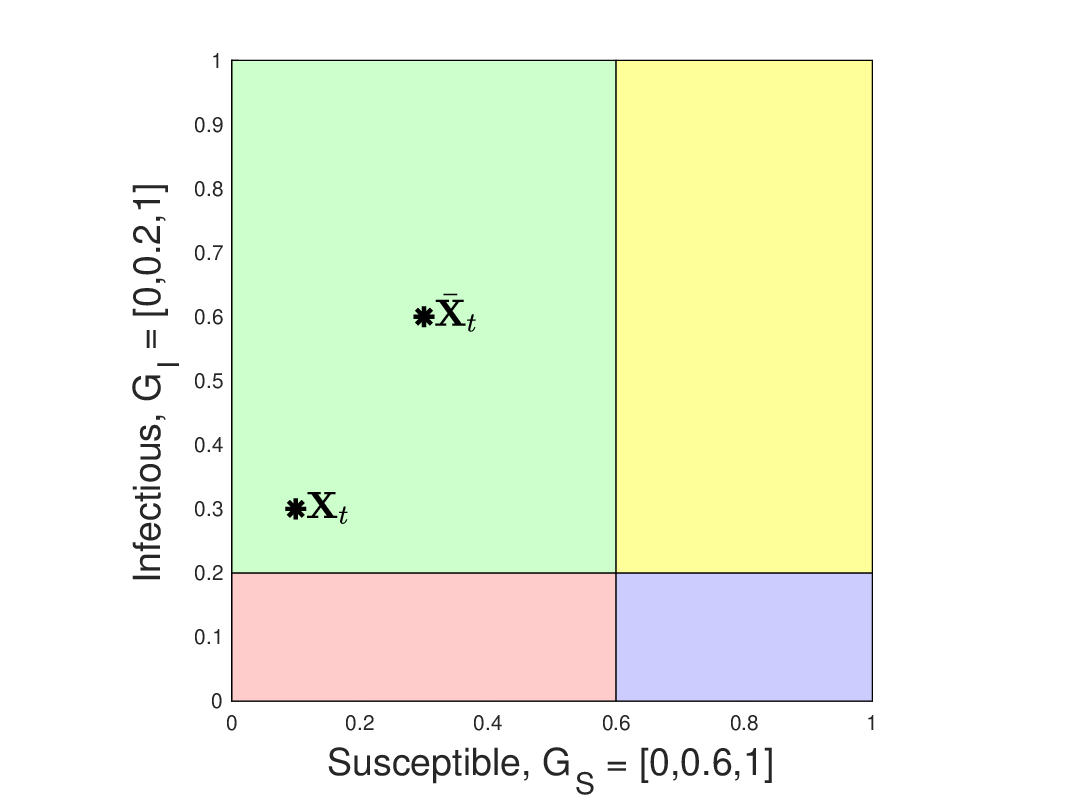}
    \caption{Four regions defined using $G=\{[0,0.6,1],[0,0.2,1]\}$ are shown in different colors. These correspond to four states:$(1):[\bar{X}_{S},\bar{X}_{I}]=[0.3,0.1];(2):[\bar{X}_{S},\bar{X}_{I}]=[0.3,0.6];(3):[\bar{X}_{S},\bar{X}_{I}]=[0.8,0.1];(4):[\bar{X}_{S},\bar{X}_{I}]=[0.8,0.6]$. For example, $\mathbf{X}_t = [0.1,0.3]$, the corresponding discretized state representation is $\bar{\mathbf{X}}_t = [0.3,0.6]$.
}
    \label{fig:discretization_vector}
\end{figure}

To discretize the continuous state space, we partition the state space $\mathcal{X}$ into a discrete set of states $\bar{\mathcal{X}}$. For each component $d$ in $\mathcal{X}$, we use the discretization vector $G_d$ to describe how the continuous state space is partitioned into discrete states. The discretization vector $G_i$ contains the maximal and minimal values of the discretized regions for component $d$. We use $G$ to represent the list of discretization vectors for all components in $\mathcal{X}$. For instance, for an SI model, if $G=\{[0,0.6,1], [0,0.2,1]\}$, we mean that group 1 (the susceptible proportion of the population) is partitioned into two regions $[0,0.6)$ and $[0.6,1]$, and the second group (infected proportion) is being partitioned into two regions $[0,0.2)$ and $[0.2,1]$. In this case, we have a total of $2\times 2 = 4$ regions. These four regions are given by $(1):X_{S}\in [0,0.6), X_{I}\in [0,0.2); (2):X_{S}\in [0,0.6), X_{I}\in [0.2,1];(3):X_{S}\in [0.6,1], X_{I}\in [0,0.2); (4):X_{S}\in [0.6,1], X_{I}\in [0.2,1]$ (shown in Figure~\ref{fig:discretization_vector}).

From these regions, we capture the discretized state space in matrix $\bar{\mathcal{X}}$, which is comprised of the Euclidean centroids of each region. The dimension of $\bar{\mathcal{X}}$ is $B\times n$ where $B$ is the number of discretization regions and $n$ is the number of components. Thus, in the example above, we would have four states. $(1):[\bar{X}_{S},\bar{X}_{I}]=[0.3,0.1]; (2):[\bar{X}_{S},\bar{X}_{I}]=[0.3,0.6]; (3):[\bar{X}_{S},\bar{X}_{I}]=[0.8,0.1]; (4):[\bar{X}_{S},\bar{X}_{I}]=[0.8,0.6]$. In this case, $\bar{\mathcal{X}}=\begin{bmatrix}
   0.3&0.1\\0.3&0.6\\0.8&0.1\\0.8&0.6
\end{bmatrix}$. Similarly, we define $\bar{\mathbf{X}}_t\in \bar{\mathcal{X}}$ to be the discretized state at time $t$ and $\{\bar{\mathbf{X}}_t\}=(\bar{\mathbf{X}}_0,...,\bar{\mathbf{X}}_N)$ to be the trajectory for the discretized state. 

With this new discretized state space, we can now define $\bar{f}(\bar{\mathbf{X}}_t,\pi_t,G)$, the disease dynamics on the discretized state space. Even though the true disease dynamics might be non-linear, we approximate the transitions on the discretized state space using a linear transition matrix. This is a reasonably good approximation if the length of $t$ is sufficiently small.

We denote this transition probability matrix as $P(\pi_t)$ for $\pi_t \in A$. $P(\pi_t)$ has the dimension of $|\bar{\mathcal{X}}|\times|\bar{\mathcal{X}}|$ where $|\bar{\mathcal{X}}|$ is the size of the state space. Then the probability of the system being in a state at time $t+1$,  $\bar{\mathbf{X}}_{t+1}$, given it was in state $\bar{\mathbf{X}}_t$ at time $t$ and policy intervention $\pi_t\in A$ is denoted as $P(\bar{\mathbf{X}}_{t+1}|\bar{\mathbf{X}}_t,\pi_t)$. It is important to note that the function $f(\mathbf{X}_t,\pi_t)$ is deterministic, whereas $P(\pi_t)$ is a stochastic matrix. The stochastic nature arises from the transition process within the discretized space. Initially, given the region of the initial state, the first transition is determined with 100\% certainty, allowing us to precisely identify the subsequent region. However, after the first transition, within any given discretized region, predicting the next region becomes uncertain as the algorithm only records past regions, not the exact points. Consequently, a transition probability matrix is employed to approximate the transition function, thereby facilitating its application to discretized MDPs.

Let $V_t(\bar{\mathbf{X}}_t)$ denote the optimal value function of the discretized state $\bar{\mathbf{X}}_t \in \bar{\mathcal{X}}, t\in T$ for the discretized infectious disease control problem. At optimality, the following must hold:
\begin{align*}
 V_t(\bar{\mathbf{X}}_t) = & \max_{\pi_t\in A} \big\{ r(\bar{\mathbf{X}}_t,\pi_t)+\lambda \sum_{\bar{\mathbf{X}}_{t+1} \in \bar{\mathcal{X}}}P(\bar{\mathbf{X}}_{t+1}|\bar{\mathbf{X}}_{t},\pi_t)V_t(\bar{\mathbf{X}}_{t+1}) \big\}
\end{align*}

\subsection{State Space Discretization Problem}
With the original system $f(\mathbf{X}_t,\pi_t)$ and state space $\mathcal{X}$, we aim to find the discretized state space $\bar{\mathcal{X}}$ and the transition matrices $P$ that approximate well the original system in that it gives a similar objective value $V_t( \bar{\mathbf{X}}_t)$, trajectories $\{\bar{\mathbf{X}}_t\}$ given $\{\pi_0,...,\pi_{N-1}\}$, and a small optimality gap. In order to do this, we need to find a suitable $G$ and map from $\bar{f}(\bar{\mathbf{X}}_t,\pi_t,G)$ to $P$.

We focus on approximating the original system by establishing an appropriate discretization approach. %An effective discretization method should consist of a small number of discretized states that consider intervention effects. 
For a discretization method to be effective, it should provide accurate estimates for regions of the state space that are more frequently visited, where these frequencies are influenced by which interventions are implemented. To do this efficiently, the discretized states should be capable of providing higher precision in areas where the state space is more likely to be visited. This can lead to a better approximation of the true disease dynamics and can thus result in a more accurate MDP solution. 

Given the function $f(\mathbf{X}_t,\pi_t)$, the initial state, the time horizon, and a sequence of policies $\{\pi_0,...,\pi_{N-1}\}$, we can calculate a trajectory $\{\mathbf{X}_t\}$. Subsequently, we require a state discretization $G$ that ensures the discretized trajectory $\{\bar{\mathbf{X}}_t\}$ closely approximates $\{\mathbf{X}_t\}$ for various initial states and policies. Therefore, our objective is to minimize the distance between the true trajectory and trajectory from the discretized model over all samples $\theta=(\mathbf{X}_0,\{\pi_0,...,\pi_{N-1}\})\in\Theta$, all policy intervention scenarios $\pi_t$, and all time, i.e., minimizing $\sum_{\theta\in \Theta}\sum_{t=1}^N  ||\mathbf{X}_t-\Bar{\mathbf{X}}_t||_2 \mid \theta$. Given a sequence of policy intervention $\{\pi_0,...,\pi_{N-1}\}$ and an initial state $\mathbf{X}_0$, we compute the true trajectory using $f(\mathbf{X}_t,\pi_t)$. We use $\bar{f}(\bar{\mathbf{X}}_t,\pi_t,G)$ to compute the trajectory from the discretization space matrix $\bar{\mathcal{X}}$.

We then map the transition function for discretized states $\bar{f}(\bar{\mathbf{X}}_t,\pi_t,G)$ to transition probability matrix $P$. Various existing techniques help to construct transition probabilities given function $\bar{f}(\bar{\mathbf{X}}_t,\pi_t,G)$. We discuss how to find a generalizable and efficient way of computing transition probabilities from $\bar{f}(\bar{\mathbf{X}}_t,\pi_t,G)$ given the state discretization in the next section. A complete list of notation used in this section is provided in Table~\ref{tab:notations}.

\section{Algorithms}\label{algo}
In this section, we present a generalizable framework for discretizing a continuous state space for use in MDP frameworks and correspondingly constructing transition probability matrices.

\subsection{Greedy Algorithm for Finding Discretization Regions (GreedyCut)}
%Due to its non-convex nature, solving the optimization problem described by equations (3)-(9) presents a significant challenge \citep{burer2012non}. 
The main objective of discretization is to design an effective approach for approximating the disease dynamics with a high level of accuracy, making such problems tractable for conventional discrete-space MDP frameworks. However, it would not be advantageous if the process of finding discretization regions itself becomes excessively costly. Therefore, our motivation is to identify a low-cost method that can produce discretization regions capable of representing the disease dynamics effectively. In particular, we are interested in outperforming a uniform discretization, which can be considered a general default discretization appropriate across many domains. 

We assume there is a budget $B$ that represents the total number of discretization regions we can have in realize of computational considerations. We use simulated initial states and policy interventions $\theta = (\mathbf{X}_0,\{\pi_0,...,\pi_{N-1}\})\in \Theta$ to find the discretization regions.

The greedy approach has been widely applied to various optimization tasks, which is easy to implement and effective at finding solutions \citep{blanchard2014greedy,wu2018greedy,zhao2021iterated}. We now propose Algorithm~\ref{euclid} (GreedyCut), a greedy-based iterative approach to finding a good discretization. A list of additional notation is listed in Table~\ref{tab:notations_algorithm}.

\begin{table}[ht]
    \centering
    \begin{tabular}{|c|c|}
    \hline
    $\Theta$ & The set of all samples used in GreedyCut \\ \hline
    $\hat{\mathbf{X}}_0$ & Samples generated in sample average approximation for generating transition probabilities\\ \hline
    $B$ & Total number of discretization regions generated\\ \hline
    $K$ & Computational costs for calculating the observed trajectory 
    given a sample\\ \hline
    $\bar{K}$ & Computational costs for calculating the discretized trajectory given a sample\\ \hline
    $\hat{K}$ & Computational costs for computing the discretized state in the next transitions in Algorithm 2 \\ \hline
    \end{tabular}
    \caption{Table of additional notation used in Section~\ref{algo}.}
    \label{tab:notations_algorithm}
\end{table}

\begin{algorithm}
\caption{Iterative Discretization for Disease Control Problems}\label{euclid}
\begin{algorithmic}[1]
\Procedure{Cost}{$\{\mathbf{X}_t\},\{\bar{\mathbf{X}_t}\}$}\Comment{$\{\mathbf{X_t}\}$ is the true trajectory, $\{\bar{\mathbf{X}}_t\}$ is the trajectory from the disretization}
\State \Return {$\sum_{t=1}^N||\bar{\mathbf{X}}_t-\mathbf{X}_t||_2^2$}
\EndProcedure
\Procedure{Cut}{$d$, $i$, $G$}\Comment{$d$ is the component we want to cut, and we want to cut the $i$-th interval in half, $i\in1,2,\ldots,G_d -1$}
\State $G_d = [0,..., G_{d,i}, (G_{d,i}+G_{d,i+1})/2, G_{d,i+1},..., 1]$ 
\State \Return{$G$}
\EndProcedure
\Procedure{Greedy}{$B$, $G$, $f(\mathbf{X}_t,\pi_t)$, $\bar{f}(\bar{\mathbf{X}}_t,\pi_t,G)$, $\theta$}\Comment{$B$ is the budget, $f(\mathbf{X}_t,\pi_t)$ is the compartmental model dynamics, $\bar{f}(\bar{\mathbf{X}}_t,\pi_t,G)$ calculates the trajectory using discretized states $\bar{\mathbf{X}}_t$, $\Theta$ is the pre-generated samples of initial states and policies, for each $\theta\in \Theta$, $\theta=(\mathbf{X}_0, \{\pi_0,...,\pi_{N-1}\})$}
\State iter\_per\_sample $= B/|\Theta|$
\For {$\theta\in\Theta$}
\For {iterations $=1:$ iter\_per\_sample}
  \State best cost  = $\infty$
  \State worst cost = $-\infty$
  \For{Component $d$}
  \For{region $i\in G_d$ }
  \State Compute $\{\mathbf{X}_t\}$ using $\mathbf{X}_t=f(\mathbf{X}_{t-1},\pi_{t-1})$
  \State Compute $\{\bar{\mathbf{X}}_t\}$ using $\bar{\mathbf{X}}_t=\bar{f}(\bar{\mathbf{X}}_{t-1},\pi_{t-1},\text{Cut($d$, $i$, $G$)})$
  \State tmp cost = Cost($\{\mathbf{X}_t\}$,$\{\bar{\mathbf{X}}_t\}$)
  \If {tmp cost $<$ best cost}
  \State best cost = tmp cost
  \State $d^* = d$
  \State $i^*=i$
  \EndIf
  \If {tmp cost $>$ worst cost}
  \State worst cost = tmp cost
  \EndIf
  \EndFor
  \EndFor
  \If{worst cost = best cost}
  \State draw a point $X_{dt}$ from $\{\mathbf{X}_t\}$
  \State Update $(d^*,i^*)$ such that $(d^*,i^*)$ satisfies $G_{d^*,i^*}\leq X_{dt}\leq G_{d^*,i^*+1} $
  \State update $G$=Cut($d^*$, $i^*$, $G$)
  \Else 
  \State update $G$=Cut($d^*$, $i^*$, $G$)
  \EndIf
  \EndFor
\EndFor
\State \Return{$G$}
\EndProcedure
\end{algorithmic}
\end{algorithm}

In Algorithm~\ref{euclid}, we have three functions. The cost function computes the sum of squared error between the trajectory from the discretized state space $\{\bar{\mathbf{X}}_t\}$ and true trajectory $\{\mathbf{X}_t\}$ from $f(X_t, \pi_t)$. Algorithm~\ref{euclid} can be adapted to any disease model without modifications to the algorithm. We compute the discretized trajectory $\{\bar{\mathbf{X}}_t\}$ using $\bar{f}(\bar{\mathbf{X}_t}, \pi_t,G)$, where the $d$-th component $\Bar{\mathbf{X}}_{dt} = \sum_{j=0}^{|G_d|-1}\mathbf{1}_{G_{d,i}\leq f(\Bar{\mathbf{X}}_{t-1},\pi_{t-1})_d < G_{d,i+1}}\frac{G_{d,i}+G_{d,i+1}}{2}$ takes the average value of the region it belongs to after the discretization. The cost function can also be customized (e.g., introduce another penalty term to emphasize certain disease compartments). 

\begin{figure}[ht]
 \centering
 \begin{subfigure}{.3\textwidth}
 \centering \includegraphics[width=2in]{G.eps}
 \caption{}
 \label{fig:discretization_G_2_a} 
 \end{subfigure} %
 \begin{subfigure}{.3\textwidth}
 \centering
 \includegraphics[width=2in]{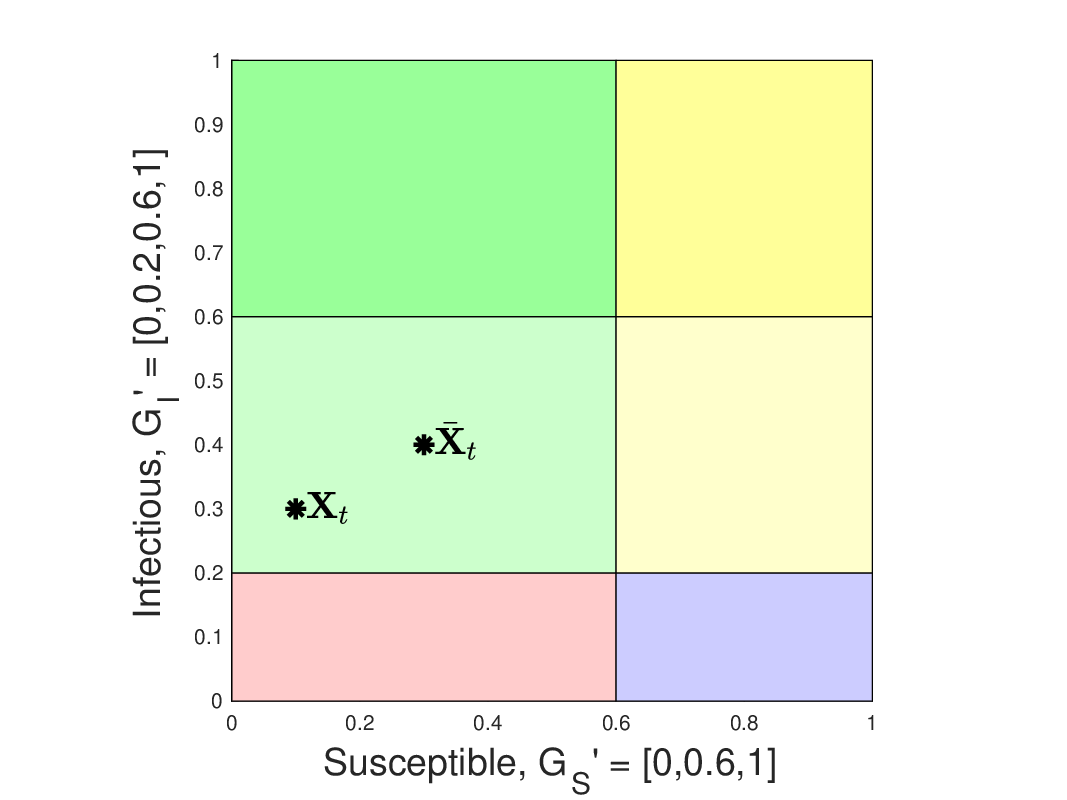}
 \caption{}
 \label{fig:discretization_G_2_b} 
 \end{subfigure} %
 \begin{subfigure}{.3\textwidth}
 \centering
 \includegraphics[width=2in]{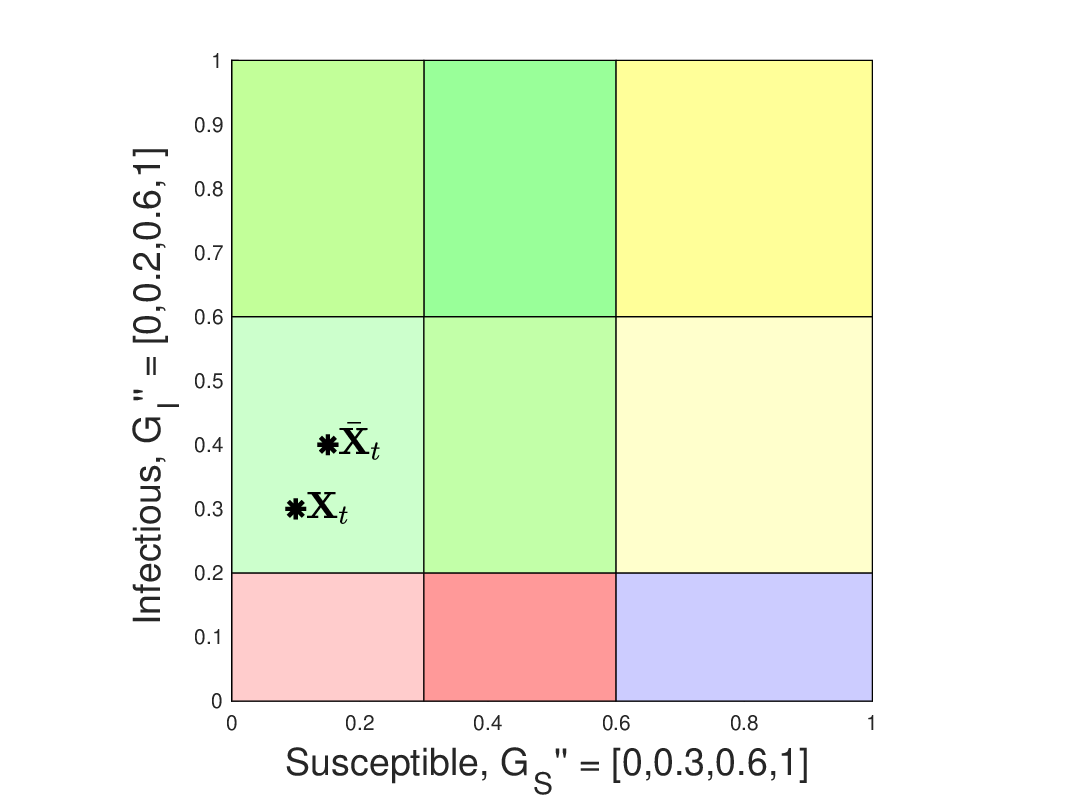}
 \caption{}
 \label{fig:discretization_G_2_c} 
 \end{subfigure} %
 \caption{Apply Cut($2$,$2$,$G$) where $G=\{[0,0.6,1], [0,0.2,1]\}$ gives new discretization regions $G'=\{[0,0.6,1], [0,0.2,0.6,1]\}$. Then apply Cut($1$,$1$,$G'$) gives new discretization regions $G''=\{[0,0.3,0.6,1], [0,0.2,0.6,1]\}$. In both $G'$ and $G''$, the $\bar{\mathbf{X}}_t$ is changed as the Euclidean centroid where $\mathbf{X}_t$ belongs to has changed. For $G$, $||\mathbf{X}_t-\bar{\mathbf{X}}_t||^2=0.11$. For $G'$, $||\mathbf{X}_t-\bar{\mathbf{X}}_t||^2=0.05$. For $G''$, $||\mathbf{X}_t-\bar{\mathbf{X}}_t||^2=0.0125$. }
 \label{fig:discretization_G_2} 
\end{figure}

The cut function divides the $i$-th region of the $d$-th component into two discretized regions, transforming a continuous range into discrete segments. For instance, consider running a single iteration of GreedyCut on $G = {[0,0.6,1], [0,0.2,1]}$ with the point $\mathbf{X}_t = (0.1,0.3)$, as illustrated in Figure~\ref{fig:discretization_G_2_a}. The initial cost is computed as $(0.3 - 0.1)^2 + (0.6 - 0.3)^2 = 0.11$. Four potential cuts are considered: $Cut(1,1,G)$ with a cost of 0.0925, $Cut(1,2,G)$ with a cost of 0.11, $Cut(2,1,G)$ with a cost of 0.11, and $Cut(2,2,G)$ with a cost of 0.05. Among these, the optimal cut is $Cut(2,2,G)$. This operation results in a new discretization, $G'$, defined as $G' = \{[0, 0.6, 1], [0, 0.2, 0.6, 1]\}$. Figure~\ref{fig:discretization_G_2_b} displays this updated discretization. After the cut, there are now six discretized regions. Similarily, running another iteration of GreedyCut on $G'$ updates the discretization to $G'' = {[0,0.3,0.6,1], [0,0.2,0.6,1]}$, with the optimal cut being $Cut(1,1,G')$ (Shown in Figure~\ref{fig:discretization_G_2_c}). After the cut, nine discretized are generated with $G''$. Consequently, the objective value decreases from 0.11 to 0.05 and then to 0.0125.

The greedy function then iteratively computes the cost of cutting one continuous range into two equal discretization regions along each component (dimension) and finds the best cut. We use the `worst cost' and `best cost' to monitor the least and most favorable cuts. Although a greedy algorithm typically does not need the `worst cost', we use it to identify cases that all cuts have equivalent costs. If each cut has the same cost, a point ($X_{dt}\in \mathbf{X}_t$) from the sampled trajectory ($\{\mathbf{X}_t\}|\theta$) will be randomly drawn, and the region that this point belongs to (component $d$ of the region $i$ of $G$ such that $G_{d,i}\leq X_{dt}\leq G_{d,i+1}$) will be cut into halves. When every cut incurs the same cost, we want to cut based on the data obtained through sampling. In general, it is unlikely that the costs for all cuts will be exactly the same; this might occur at the beginning of the algorithm when each discretized state encompasses a large range and the approximation will not improve if cut only once. An alternative solution is to generate a few random cuts at the beginning and then run GreedyCut without recording the `worst cost'. Through this process, in total, $|G|=\sum_d |G_d|=B$ discretization regions will be generated. In Section~\ref{algo_evaluations}, we will assess the effectiveness of this algorithm by comparing the sum of squared errors between the Markovian trajectory (the expected trajectory computed using probabilities and discretization regions) and the observed trajectory with different discretization methods. This comparison will help us evaluate how well GreedyCut estimates the actual disease trajectory using transition probability matrices.
\subsubsection{Complexity Analysis}\label{comp_algo1}
In Algorithm~\ref{euclid}, if we assume that computing $\{{\mathbf{X}_t}\}$ and $\{\bar{\mathbf{X}}_t\}$ using $f(\mathbf{X}_t, \pi_t)$ and $\bar{f}(\bar{\mathbf{X}}_t, \pi_t, G)$ given one sample requires $K$ and $\bar{K}$ operations respectively, we can analyze the total number of operations performed by the GreedyCut algorithm. Since the computational costs for generating new discretization regions,  comparing costs, and computing costs are relatively small compared to computing $\{\mathbf{X}_t\}$ and $\{\bar{\mathbf{X}}_t\}$ in our problem, we assume these costs are negligible compared with other costs. At each iteration, the GreedyCut algorithm enumerates through each of $b$ regions including all components $d$, where $1 \leq b \leq B$. For each discretization, the algorithm performs computations of $\{\mathbf{X}_t\}$ and $\{\bar{\mathbf{X}}_t\}$, which have time complexities of $K$ and $\bar{K}$ operations respectively. The total number of iterations depends on $B$ and $|\Theta|$. It equals to $|\Theta|\times \frac{B}{|\Theta|}=B.$

The total number of operations for the GreedyCut algorithm can be estimated as the sum of operations over all iterations. This can be expressed as $\sum_{b=1}^{B} b(K+\bar{K})=\frac{B(B+1)}{2}(K+\bar{K})$. Therefore, the complexity of the GreedyCut algorithm is on the order of $\mathcal{O}(B^2(K + \bar{K}))$. This implies that the complexity grows exponentially with the budget of discretization regions,  given by $B$. The time complexity increases quadratically with $B$ while linearly with the number of operations required for each region, represented by $(K + \bar{K})$.

As a result, the computational complexity of the algorithm grows rapidly as the number of discretization regions increases. This highlights the exponential relationship between the complexity and the desired level of granularity in the discretization process. 

\subsection{Constructing a Corresponding Transition Matrix}

Once a suitable discretization of a continuous state space has been constructed, we additionally need a transition matrix between these discretized states to capture the dynamics for use in an MDP framework. To do this, we use a well-known sample average approximation (SAA)-based algorithm (hereafter referred to as Algorithm 2), commonly utilized in various algorithms \citep{kim2015guide,bertsimas2018robust,zhang2024sampling}. In our approach, SAA draws samples from each discretized state to determine the frequency of transitions to subsequent states via the function $f(\mathbf{X}_t,\pi_t)$ and $G$ \citep{hernandez2012discrete}. The detailed algorithm steps are provided in the Appendix.

In Algorithm 2, we draw $c$ samples within each region in $G$ and count the frequency of the transitions from the current region to other regions using policy $\pi_t$ and $f(\mathbf{X}_t,\pi_t)$. By sampling and counting the transitions from the original system, we approximate the underlying transition probabilities directly. Creating approximated transition probability matrices in this way offers a practical approach to capturing the essential dynamics of the system and enables efficient decision-making at the population level.

\subsubsection{Complexity Analysis}
In Algorithm 2, we generate $c$ samples $B|A|$ times, where $B$ denotes the number of discretization regions and $|A|$ represents the size of the action space. Assuming that computing and locating each transition in the appropriate region take $\hat{K}$ operations, we can analyze the time complexity of Algorithm 2. The number of operations performed by the algorithm then is $\mathcal{O}(c\hat{K}B|A|)$.

Furthermore, in most disease control problems, such as COVID-19 mitigation strategies like lockdown, social distancing, and face masks, the size of the action space $|A|$ is typically small. This implies that the algorithm's time complexity is primarily influenced by the number of samples $c$, the number of discretization regions $B$, and the operations $\hat{K}$ needed for computing and locating transitions. 

\section{Numerical Examples}\label{numeric_discretization}
In this section, we first showcase our proposed framework for reformulating a SIR model to an infectious disease control MDP framework for supporting public health decisions around social distancing policy. We then demonstrate the utility of this framework with an example of COVID-19 in Los Angeles County, drawing from empirical data of case counts in 2020. 

We benchmark the outcome of our method (we refer to the `GreedyCut discretization method' hereafter) in both examples by comparing our model outcomes to those of several other frameworks: a uniform discretization framework, a discretization framework based on expert opinion (we refer to the `expert discretization method' hereafter), and an InverseProportional discretization framework motivated by \cite{zhou2010solving}. In the uniform discretization framework, we discretize the entire state space uniformly using the same number of discretization regions as used in the GreedyCut discretization method. In the Expert discretization framework, the state spaces for susceptible and recovered individuals are discretized uniformly. However, the proportion of the population that is infected is only uniformly discretized between $[0,0.4]$ instead of $[0,1]$. This is based on the observation that the proportion of the population infected at the same time does not exceed 40\% in most infectious diseases \citep{biggerstaff2014estimates,pei2021burden}. In the InverseProportional discretization framework, the state space is discretized by assigning higher resolution to the more frequently visited states based on the simulation results. The transition probabilities for all methods are generated using Algorithm 2. In the second example, we additionally compare our model outcomes to the empirical status-quo policy in Los Angeles in 2020 to demonstrate the improvement our method can achieve.

\subsection{Example 1: A Simple SIR Model}\label{sir}
The SIR model tracks the proportion of the population that is susceptible ($S$), infected ($I$), and recovered ($R$) at each time $t$. We use a discrete time model where the SIR model can be described using a system of difference equations \citep{allen1994some}:
 \begin{align*}      S_{t+1} & =S_t -\beta  S_t I_t \\
     I_{t+1}&=I_t +\beta  S_t I_t -\gamma I_t \\
     R_{t+1}&=R_t+\gamma I_t 
 \end{align*}

The parameter $\beta$ is the rate at which disease transmits from the infected to susceptible population proportions, and is dependent on the average contact rate and probability of transmission given a discordant contact. Similarly, $\gamma$ is the recovery rate.

Typically, at the beginning of an epidemic, the exact proportion of the population that is infected may be unknown. We use $\mathbf{X}_0=[X_{S0},X_{I0},X_{R0}]=[S_0,I_0,R_0]$ to denote the initial state at the first decision epoch. We assume that while the exact initial state is unknown, we know an upper and lower bound on each of the compartments. We use $\underline{S},\bar{S}$, $\underline{I},\bar{I}$, and $\underline{R},\bar{R}$ to denote the upper and lower bound on initial states $S_0$, $I_0$, and $R_0$, respectively.

Suppose at each time $t$, the health department can choose to implement a social distancing policy (a ``lockdown") until time epoch $t+1$ that reduces the transmission rate $\beta$. We assume there are a finite number of periods $N$. 

The decision maker wishes to minimize the negative health outcomes and economic and social costs of implementing a lockdown policy. To capture this objective, at each decision epoch, we let the cost be $r(\mathbf{X}_t,\pi_t) = I_t + u(\pi_t)$, the proportion of the population infected ($I_t$) plus some time-invariant dis-utility value $u(\pi_t)$ that captures the economic and social costs that are only incurred when the intervention is in effect and zero otherwise.

 Throughout this section, we refer to this discrete time system as the ground-truth system, and we will construct our discretized MDP framework based on this. We assume no discounting in the objective ($\lambda$=1). The objective of the MDP is therefore to minimize the total costs over the whole time horizon.

\subsubsection{Inputs}\label{example_inputs}
To evaluate this example, we let the transmission rate ($\beta$) be 1.4 and the recovery rate ($\gamma$) be 0.49. The decision interval ($\Delta t$) is a week, and the time horizon ($N$) is ten weeks. Implementing a lockdown will incur an economic and social cost, but it is unclear how this dis-utility can be quantified in reality. For simplicity, we assume the dis-utility is 0.03 if lockdown was implemented and 0 otherwise ($u(\text{lockdown})=0.03, u(\text{no lockdown})=0$). 

During the early stages of a pandemic, there is typically a large population in the susceptible category, while only a small population is infected. Therefore, we choose the initial states to be uniformly distributed within the upper and lower bounds for each compartment to be: $[\underline{S},\bar{S}]=[0.7,0.99]$, $[\underline{I},\bar{I}]=[0.01,0.1]$, and $[\underline{R},\bar{R}]=[0,0.29]$.

In the GreedyCut discretization method, for each sample ($\theta$) generated, ten iterations are run to generate ten additional discretization regions (lines 9 - 24 in Algorithm~\ref{euclid}). To generate samples $\theta$, $\mathbf{X}_0$ are generated uniformly from the region above, and $\pi$ is a vector with ten random binary variables to indicate the policy intervention (0 -- no lockdown, 1 --  lockdown). In Algorithm 2, we generate $c=1000$ samples to compute two transition probability matrices to correspond to the no lockdown and lockdown policies, respectively.  The reason for choosing $c = 1000$ is that we observed the matrix stabilizes once $c$ exceeds 1000.

\subsubsection{MDP Solutions }
We compare MDP solutions among the GreedyCut, InverseProportional, expert, and uniform discretization methods on 90, 150, 300, and 1200 discretization regions. 

To evaluate our algorithm's performance, we create 300 samples (we denote the set of all samples as $\check{\mathbf{X}}_0$) by selecting the initial susceptible proportion from 0.7 to 0.99 with stepsize of 0.01 and initialize the infected proportion from 0.001 to 0.01 with 0.001 stepsize. We chose 3,000 state-time pairs because we observed that performance does not vary significantly when the number of state-time pairs is changed. By enumerating each pair of $S$ and $I$, we can have a total of 300 different possible initial states (if $S+I>1$, we will renormalize each compartment). We compute the following metrics for both GreedyCut and uniform discretization methods on 90, 150, 300, and 1200 discretization regions :
\begin{itemize}
    \item ACC: accuracy in matching the percentage of optimal actions by comparing discretized MDP with brute force (ground-truth) solution over each state-time pair (a total of 3000 state-time pairs). ACC $= 1-\frac{\#mismatch}{3000}$.
    \item MSE: mean squared error between the optimal value of the discretized MDP ($\bar{V}^*_0$) and the brute force solution ($V^*_0$) on the first decision epoch over all states. MSE $=\mathbb{E}_{\mathbf{X}_0\in \check{\mathbf{X}}_0}[||\bar{V}^*_0(\mathbf{X}_0)-V^*_0(\mathbf{X}_0)||^2]$.
    \item E2: relative mean absolute error on the first decision epoch over all states. E2 $=\mathbb{E}_{\mathbf{X}_0\in \check{\mathbf{X}}_0}[\frac{|\bar{V}^*_0(\mathbf{X}_0)-V^*_0(\mathbf{X}_0)|}{V^*_0(\mathbf{X}_0)}]$.
    \item Opt. Gap: average of the relative difference between the optimal value of brute force solution and the value of running optimal policy from discretized MDP on the true disease model ($\tilde{V}_0$) on the first decision epoch. Opt. Gap $= \mathbb{E}_{\mathbf{X}_0\in \check{\mathbf{X}}_0}[\frac{|\tilde{V}_0(\mathbf{X}_0)-V^*_0(\mathbf{X}_0)|}{V^*_0(\mathbf{X}_0)}]$.
\end{itemize}

\begin{figure}[ht]
 \centering
 \begin{subfigure}{.49\textwidth}
 \centering \includegraphics[width=3in]{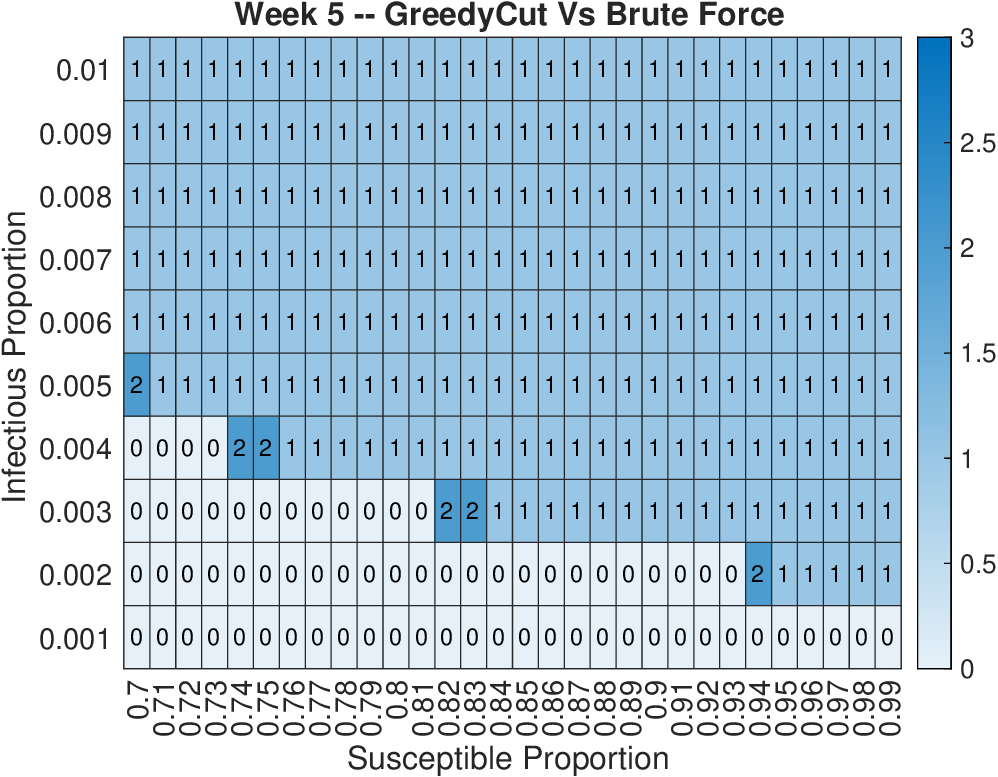}
 \caption{}
 \end{subfigure} %
 \begin{subfigure}{.49\textwidth}
 \centering
 \includegraphics[width=3in]{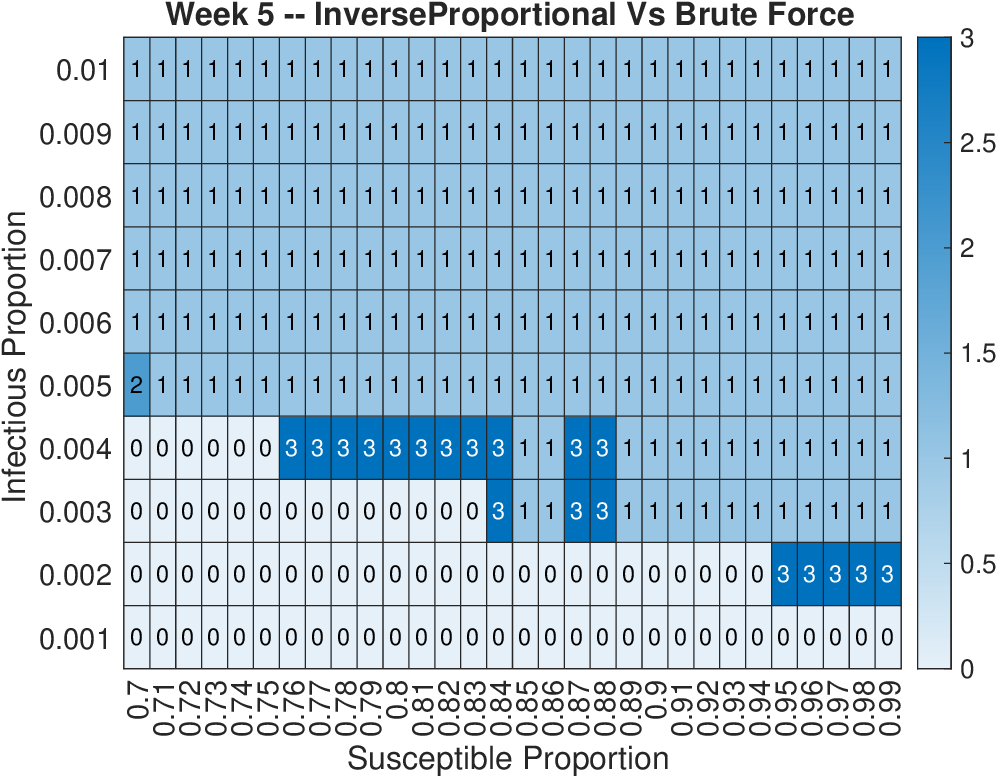}
 \caption{}
 \end{subfigure} %
   \begin{subfigure}{.49\textwidth}
 \centering
 \includegraphics[width=3in]{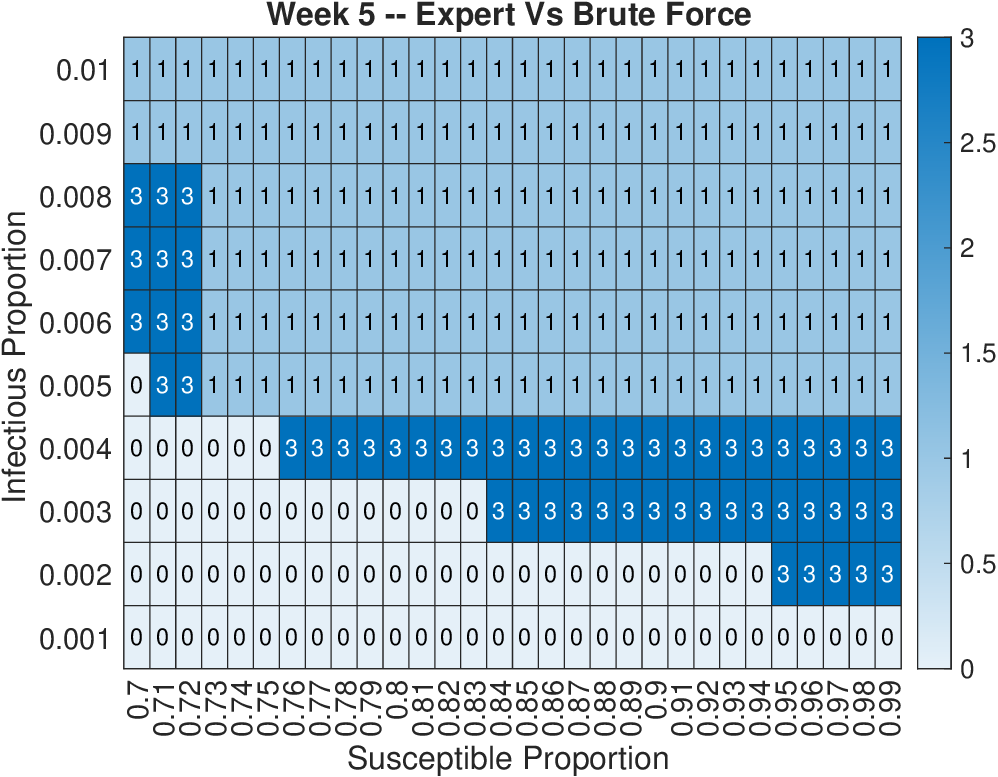}
 \caption{}
 \end{subfigure} %
  \begin{subfigure}{.49\textwidth}
 \centering
 \includegraphics[width=3in]{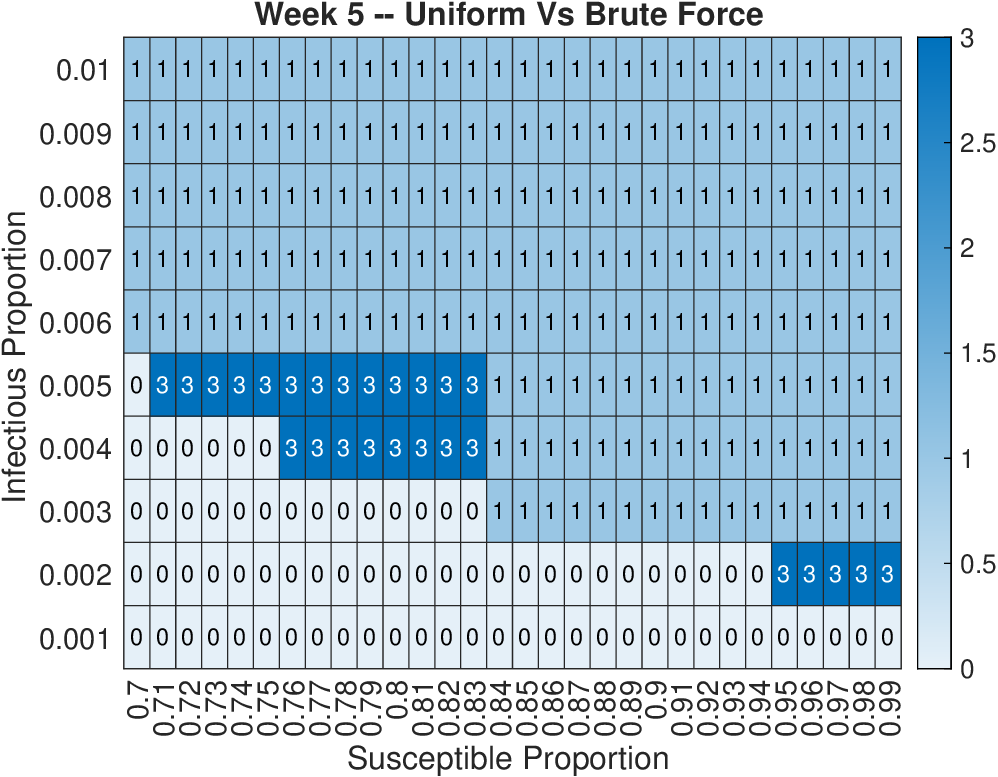}
 \caption{}
 \end{subfigure} %
 \caption{We compare the optimal solution at $t=5$ across different states using the GreedyCut and uniform discretized MDPs against the ground truth optimal solution found using brute force methods. (a): Optimal solution from the GreedyCut discretized MDP compared to the brute force solution; (b): Optimal solution from the uniform discretized MDP compared to the brute force solution. ( 0 -- both models recommend not implementing lockdown; 1 -- both models recommend implementing lockdown; 2 -- the brute force method recommends not implementing lockdown while the other method recommends lockdown [type 2 error]; 3 -- the brute force method recommends implementing lockdown while the other method recommends not implementing lockdown [type 1 error].)}
 \label{fig:compare_week5}
\end{figure}

The GreedyCut discretized MDP is able to generate solutions with a higher accuracy than the uniform discretized MDP. We compare the solutions from GreedyCut and other discretized MDPs against the brute force solution in the fifth week ($t=5$, shown in Figure~\ref{fig:compare_week5}). We chose the fifth week to better illustrate the outcomes, as the decision boundary is less illustrative in earlier and later time periods. All state pairs recommend the same optimal actions. In the fifth week, we compare optimal actions across 300 states. The results indicate that the GreedyCut discretized MDP has six mismatches, whereas the uniform discretized MDP has 26 mismatches, the expert discretized MDP has 56 mismatches, and the InverseProportional discretized MDP has 20 mismatches when compared with the brute force solution. 

Moreover, the number of mismatches may be worse with the other discretized MDPs compared to GreedyCut. In the GreedyCut discretized MDP, all six mismatches belong to the case where the brute force solution recommends not to lockdown while the GreedyCut discretized MDP recommends implementing lockdown. However, in all other discretized MDPs, most mismatches belong to the case where the brute force solution recommends lockdown while the discretized MDP recommends not implementing lockdowns. In infectious disease control, failing to implement a lockdown when it is necessary can cause a rapid increase in the proportion of infected cases. Therefore, error in this direction may be practically worse than in the converse direction. We see this illustrated in the optimality gap among four discretized MDPs (last columns in Table \ref{tab:mdp_table}), which measures the distance between solutions from discretized MDPs and the true optimal solution. Here we see that GreedyCut MDP achieves a much lower optimality gap for this reason. 

\begin{table}[ht]
 \footnotesize
 \centering
 \begin{tabular}{|c|c|c|c|c|c|c|c|c|}
 \hline
  $B$ &  \multicolumn{4}{c|}{ACC} & \multicolumn{4}{c|}{MSE}\\
  \hline
   & GreedyCut &InverseProportional&Expert&Uniform&  GreedyCut&InverseProportional&Expert & Uniform\\
   \hline
   90&0.9657&0.8687&0.8927&0.8120&4.3239e-04&0.0040&0.0084& 0.0580\\
   \hline
   150&0.9850&0.9073&0.8753&0.8820&1.7896e-04&0.0012&0.0037&0.0169\\
   \hline
   300&0.9787& 0.9187&0.8780&0.8613&6.5032e-05& 4.0714e-04&0.0011&0.0054\\
   \hline
   1200&0.9893&0.9767&0.9423&0.9150&3.0731e-06&4.8719e-05&4.9485e-04&4.3529e-04\\
   \hline
  \hline
  $B$&\multicolumn{4}{c|}{E2}& \multicolumn{4}{c|}{Opt. Gap}\\
\hline
& GreedyCut &InverseProportional&Expert& Uniform &  GreedyCut &InverseProportional&Expert&Uniform\\
 \hline
90&0.0689&0.2134&0.3115&0.8846&0.0033& 0.0194&0.0364&0.0954\\
 \hline
150&0.0435&0.1180&0.2060&0.4946&0.0011&0.0093&0.0164&0.0583\\
 \hline
300&0.0233& 0.0537&0.1012&0.2487&0.0018& 0.0071&0.0111&0.0251\\
 \hline
 1200&0.0048&0.0220&0.0775&0.0654&4.5576e-04&8.4937e-04&0.0023&0.0072\\ 
 \hline 
 \end{tabular}
 \caption{Comparison on MDP solutions}
 \label{tab:mdp_table}
\end{table}

The GreedyCut discretization method outperforms other discretization methods across different evaluation metrics across all time periods and different discretization budgets. Table~\ref{tab:mdp_table} shows the comparison among three discretization methods. Also, the GreedyCut discretization method has higher accuracy (approximately 10\% more compared to the uniform discretization method) in matching the optimal actions from the brute force (ground-truth) solution over all state-time pairs. The GreedyCut discretization method is able to generate accurate recommendations on policy interventions even with a small number of discretization regions. Additionally, this method is able to provide a closer approximation of the objective value in both MSE and E2 metrics across all discretization regions. When the number of discretization regions is small, the GreedyCut algorithm has an MSE that is under 1\% of the MSE generated using the uniform discretization approach. Similarly, under these conditions, the GreedyCut algorithm's E2 remains below 10\% of the E2 from the uniform discretization method. Moreover, the GreedyCut algorithm outperforms the uniform discretization method in reducing the optimality gap. The optimality gap ranges from 0.1\% to 0.33\% with different numbers of discretization regions in the GreedyCut algorithm, compared with its ranges from 0.72\% to 9.54\% in the uniform discretization method. We also implemented the discretization-free algorithm from \citep{li2005lazy}; however, it achieved only about 50\% accuracy, which is significantly lower than all other methods. Therefore, we excluded this method from our comparisons. This poor performance might be due to the short time horizon of the problem and the need for accurate piece-wise constant approximations of the value and transition functions.

As expected, with a small number of discretization regions, the difference in performance between the GreedyCut and the uniform discretization methods is large. The performance gap shrinks when the number of discretization regions increases, as uniform discretization regions naturally benefit from smaller discretized regions -- higher resolution. Additionally, to assess the robustness of GreedyCut against different cost functions, we implemented GreedyCut with the mean absolute error cost function. Although it performed slightly worse than the original GreedyCut, it still outperformed all benchmarking solutions across all numbers of discretization regions considered (results in the Appendix Table~\ref{tab:appendix_tab1}).  
\\
\noindent\emph{\textbf{Interpretation of Cuts Generated.}} We compare the cuts generated by GreedyCut with those from other discretization methods. We found that the GreedyCut algorithm allocated most of its cuts to states that are more likely to be visited based on sample trajectories. Unlike the InverseProportional method, which focuses solely on frequently visited states, GreedyCut also made cuts in less frequently visited but still possible regions, potentially offering a more reliable and robust discretization (the distribution of cuts can be found in Appendix Figure~\ref{fig:distribution_cuts}).

\noindent\emph{\textbf{Run Time Outcomes.}} To understand how much time is needed to construct an MDP using the GreedyCut discretization method and Algorithm 2, we compare the run time of Algorithm~\ref{euclid} and Algorithm 2 with different numbers of discretization regions using Matlab 2022b on a laptop with 16 GB memory and Apple M1 pro chip.
\begin{figure}[ht]
 \centering 
 \includegraphics[width=3.5in]{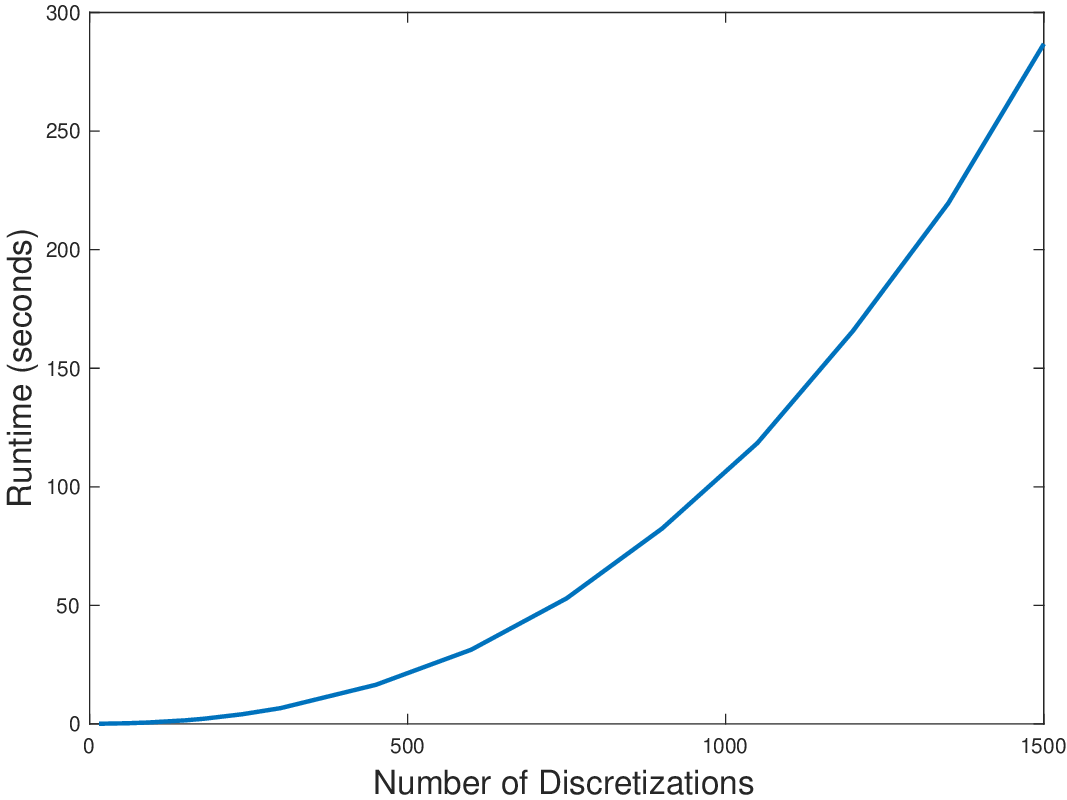} 
 \caption{Runtime of Algorithm~\ref{euclid}}
 \label{fig:runtime} 
\end{figure}

There is an exponential relationship between the number of discretization regions and the algorithm runtime (see Figure~\ref{fig:runtime}), consistent with the time complexity analysis in Section~\ref{comp_algo1}. The curvature of the exponential function will depend on the complexity of the disease model $f(\mathbf{X}_t,\pi_t)$; with more compartments or population stratifications, the total runtime may be larger for a similar number of discretization regions.

\begin{table}[ht]
 \footnotesize
 \centering
 \begin{tabular}{|c|c|}
 \hline
 $B$ & Runtime (hours):\\
 \hline
 90&0.12\\
 \hline
 150&0.34\\
 \hline
 $>$300&$>$1\\
 \hline 
 \end{tabular}
 \caption{Runtime of Algorithm 2}
\label{tab:transitions_runtime}
\end{table}

The runtime of generating transition matrices is much more costly compared with generating the discretization regions when number of discretization regions ($B$) is small for both GreedyCut and other discretization methods. Table~\ref{tab:transitions_runtime} shows the runtime of generating transitions using Algorithm 2 for the GreedyCut discretization method (other discretization methods should have the same runtimes as there are same number of iterations needed). The runtime exceeds one hour with 300 discretization regions given $c=$1000. Therefore, when $B$ is small, the total runtime of constructing an approximate MDP is roughly the time for generating transitions using Algorithm 2. In this case, the time required to construct the MDP using the GreedyCut discretization method is similar to that of using the uniform discretization method. This is because the runtime of Algorithm \ref{euclid}, which is only needed for GreedyCut and not for the uniform discretization, is negligible compared to the runtime of Algorithm \ref{sample}. The total runtime of each algorithm, including discretization, generation of the transition probability matrix, and solving the MDP can be found in the Appendix.
%In this case, the time that it takes to construct the MDP using the GreedyCut discretization method is close to that of using other discretization methods -- because the runtime of the Algorithm~\ref{euclid} (which is only needed for GreedyCut and not other discretization regions ) is negligible compared to the runtime of Algorithm 2.

\subsubsection{General Algorithmic Evaluations} \label{algo_evaluations}
In this section, we evaluate the GreedyCut discretization method's performance on generating discretization regions that approximate the disease dynamics $\{\mathbf{X}_t\}$ and Algorithm 2's performance on generating transition probability matrices to approximate the discretized trajectories $\{\bar{\mathbf{X}}_t\}$. We first examine the performance of Algorithm 2 to highlight its capability to generate precise transition probabilities. These probabilities are crucial for describing the discretized trajectory across various discretization settings. Subsequently, we assess the performance of the GreedyCut discretization method by comparing the Markovian trajectories (using transition probabilities from Algorithm 2) against the actual trajectory $\{\mathbf{X}_t\}$ between the GreedyCut algorithm and other discretization methods.
\\
\noindent\textit{\textbf{How Accurate Are the Generated Transition Probabilities?}} To evaluate the accuracy of generated transition probabilities from Algorithm 2, we draw samples and evaluate computed trajectories compared with trajectories $\{\bar{\mathbf{X}}_t\}$ from $\bar{f}(\bar{\mathbf{X}}_t,\pi_t,G)$ to eliminate the influence of the quality of the discretization algorithm.  

For evaluation, we uniformly draw $100$ samples ($\hat{\theta} \in \hat{\Theta}$) that consist of initial states within the upper and lower bound of the proportions in each compartment and a sequence of policy interventions for each initial state. For each evaluation sample $\hat{\theta}$, we compute the discretized trajectory $\{\bar{\mathbf{X}}_t\}$ using $\bar{f}(\bar{\mathbf{X}}_t,\pi_t,G)$. To obtain the Markovian trajectories from the discretized Markov model, we use an initial belief $b_0 = \mathbf{e}_i$ where all entries of $b_0$ are zero except for $i$-th entry (corresponding to $\mathbf{X}_0$) which has value one. This indicates we know 100\% the initial state of the discretized Markov model. Then we update the belief $b_t = P(\pi_t)b_{t-1}$ over time. To compute the expected proportion of people on each time $t$ (Markovian trajectory at time $t$, $\Tilde{\mathbf{X}}_t$), we use the weighted average over the belief vector at time $t$, e.g., $\tilde{\mathbf{X}}_{t}=b_t^T\bar{\mathcal{X}}$. We then compute the cost $\mathbb{E}_{\hat{\theta} \in \hat{\Theta}}[\sum_{t=1}^N||\tilde{\mathbf{X}}_{t}-\bar{\mathbf{X}}_t||^2_2]$ to evaluate how close is Algorithm 2 able to generate reliable transition probability matrices.
 
\begin{figure}[ht]
 \centering 
 \includegraphics[width=4in]{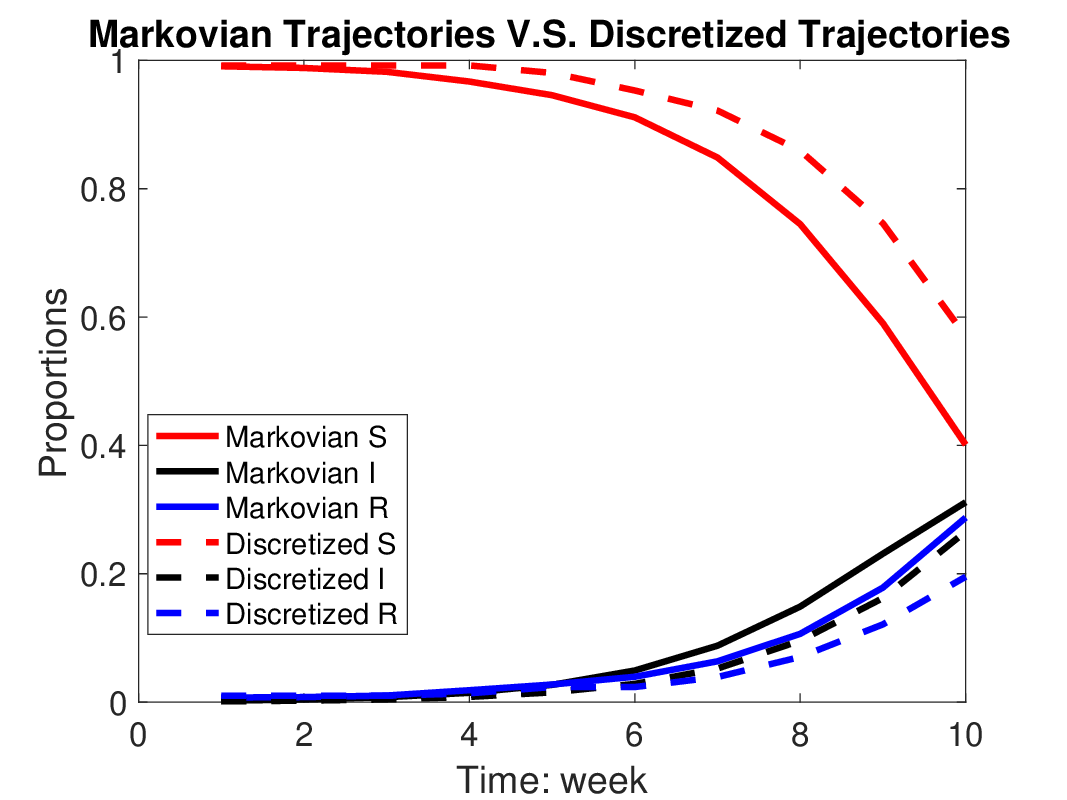} 
 \caption[Comparison between trajectories.]{Comparison between trajectories generated from Algorithm 2 given discretization regions and trajectories generated from discretized states. For each compartment S, I, and R, both trajectories are close to each other. }
 \label{fig:figure3.1} 
\end{figure}

The trajectories obtained from Algorithm 2 closely align with those generated from discretized states for each compartment. As shown in Figure~\ref{fig:figure3.1}, we compared the trajectories obtained from Algorithm 2 using 300 discretization regions with the $\bar{f}(\bar{\mathbf{X}}_t,\pi_t,G)$ trajectories generated from the same 300 discretized states. 

We observed that as the number of discretization regions increases, Algorithm 2 is capable of generating transitions that closely resemble the dynamics of the disease, represented by $\{\bar{\mathbf{X}}_t\}$. In Table~\ref{tab:transitions_mse}, we show a comparison of the sum of squared error between the Markovian trajectory and discretized trajectory ($\mathbb{E}[\sum_{t=1}^N||\tilde{\mathbf{X}}_{t}-\bar{\mathbf{X}}_t||^2_2]$) for total budget of discretization regions ($B$) of 90, 150, 300, and 1,200 using the GreedyCut discretization method. We observe that the error consistently decreases as the number of discretization regions increases.
%In Table~\ref{tab:transitions_mse}, we show a comparison of the cost $\mathbb{E}[\sum_{t=1}^N||\tilde{\mathbf{X}}_{t}-\bar{\mathbf{X}}_t||^2_2]$ for $B$ of 90, 150, 300, and 1,200 using the GreedyCut discretization method. We observe that the population proportions generated using the Markovian and discretized processes closely align, even over time.

With fewer discretization regions,  each individual discretization possesses a larger range, making it more difficult for samples drawn from these discretization regions to transition accurately between decision epochs, leading to more error in approximating $\{\bar{\mathbf{X}}_t\}$. On the other hand, when a larger number of discretization regions is employed, each discretization exhibits a smaller range. By drawing a sufficient number of samples, it becomes possible to provide a more precise description of $\{\bar{\mathbf{X}}_t\}$. These findings highlight the algorithm's reliability and accuracy in capturing the system's dynamics.

\begin{table}[ht]
 \footnotesize
 \centering
 \begin{tabular}{|c|c|}
 \hline
 $B$ & Mean Squared Error Between Discretized Trajectory and Markovian Trajectory [95\% uncertainty interval] \\
 \hline
 90&$0.1598 \quad [0.1411,0.1785]$\\
 \hline
 150&$0.1173 \quad [0.1020,0.1326]$\\
 \hline
 300&$0.1086 \quad [0.0940,0.1232]$\\
 \hline
 1200&$0.1067 \quad [0.0923,0.1211]$\\ 
 \hline 
 \end{tabular}
 \caption{Mean squared error for the trajectories given different numbers of discretization regions }
 \label{tab:transitions_mse}
\end{table}

\noindent\textit{\textbf{How Accurate Are the Discretization Regions Generated from the GreedyCut Discretization Method? }} Next, to evaluate the quality of discretization regions generated from the GreedyCut discretization method, we draw samples and calculate the mean squared error across all samples. This error is measured between the actual trajectory ($\{\mathbf{X}_t\}$) and anticipated Markovian trajectory ($\{\tilde{\mathbf{X}}_t\}$), using the transition matrix created in Algorithm 2 using the discretization regions generated from Algorithm~\ref{euclid} based on the same 100 samples for evaluating Algorithm 2.
\begin{figure}[ht]
 \centering
 \begin{subfigure}{.49\textwidth}
 \centering \includegraphics[width=2.5in]{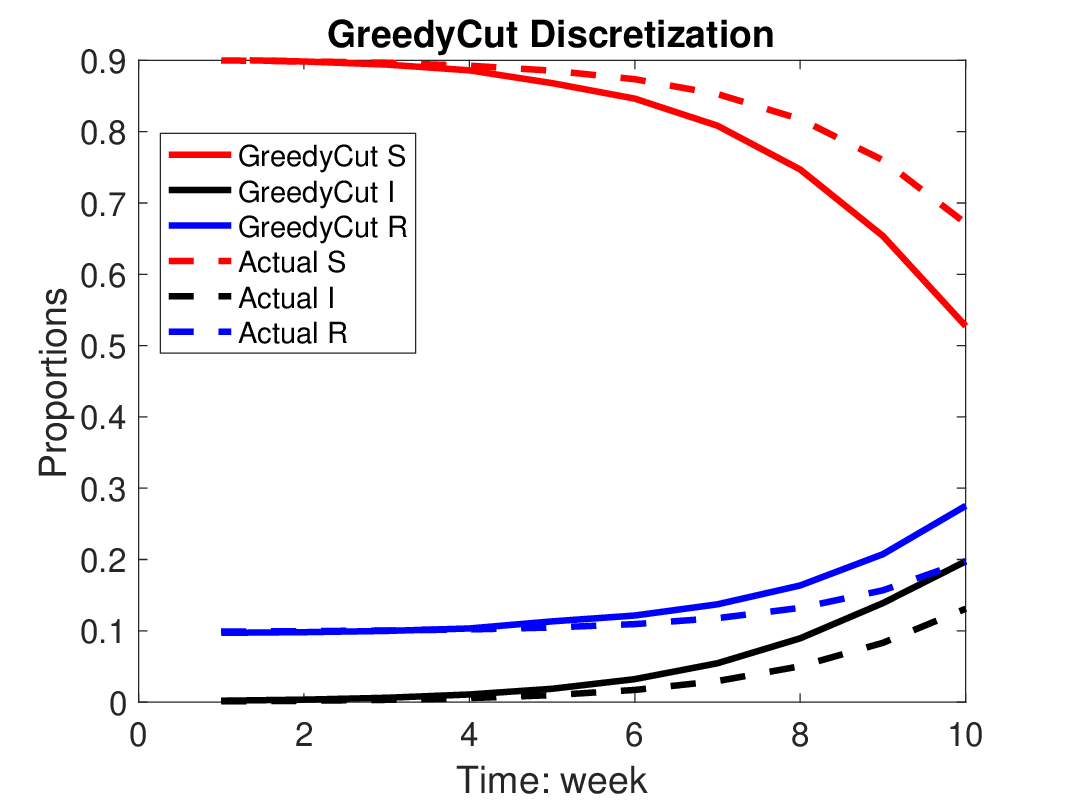}
 \caption{}
 \end{subfigure} %
 \begin{subfigure}{.49\textwidth}
 \centering \includegraphics[width=2.5in]{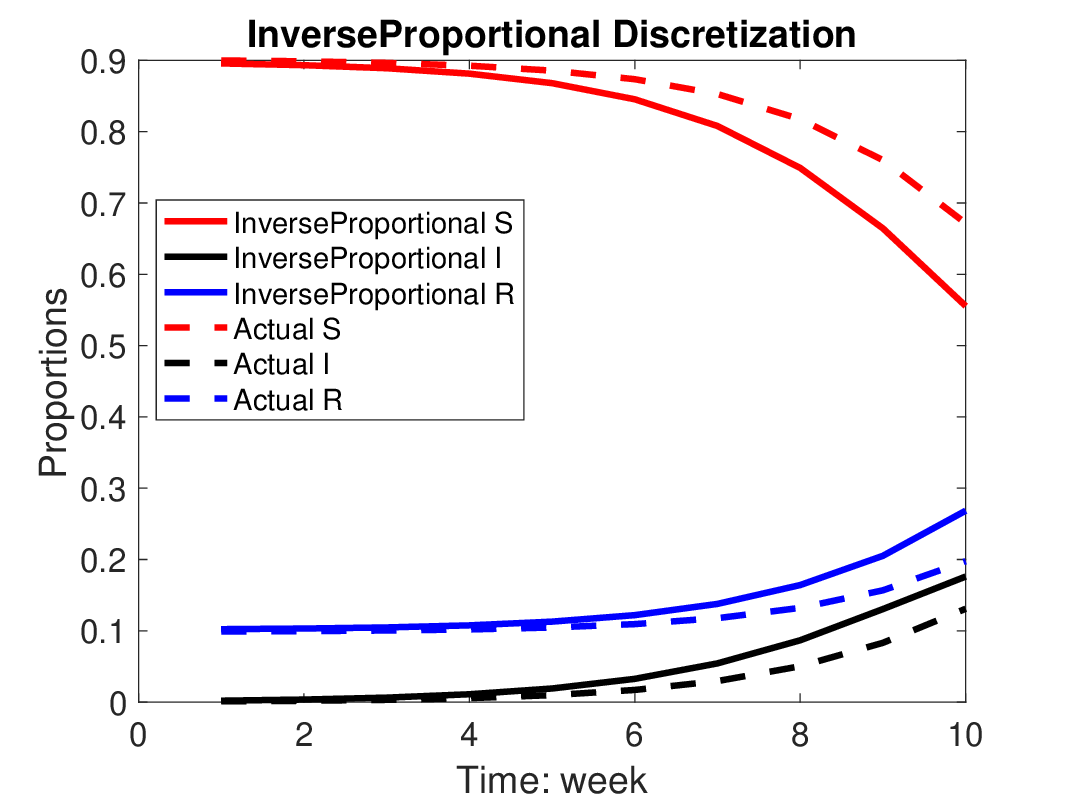}
 \caption{}
 \end{subfigure} %
 \begin{subfigure}{.49\textwidth}
 \centering
 \includegraphics[width=2.5in]{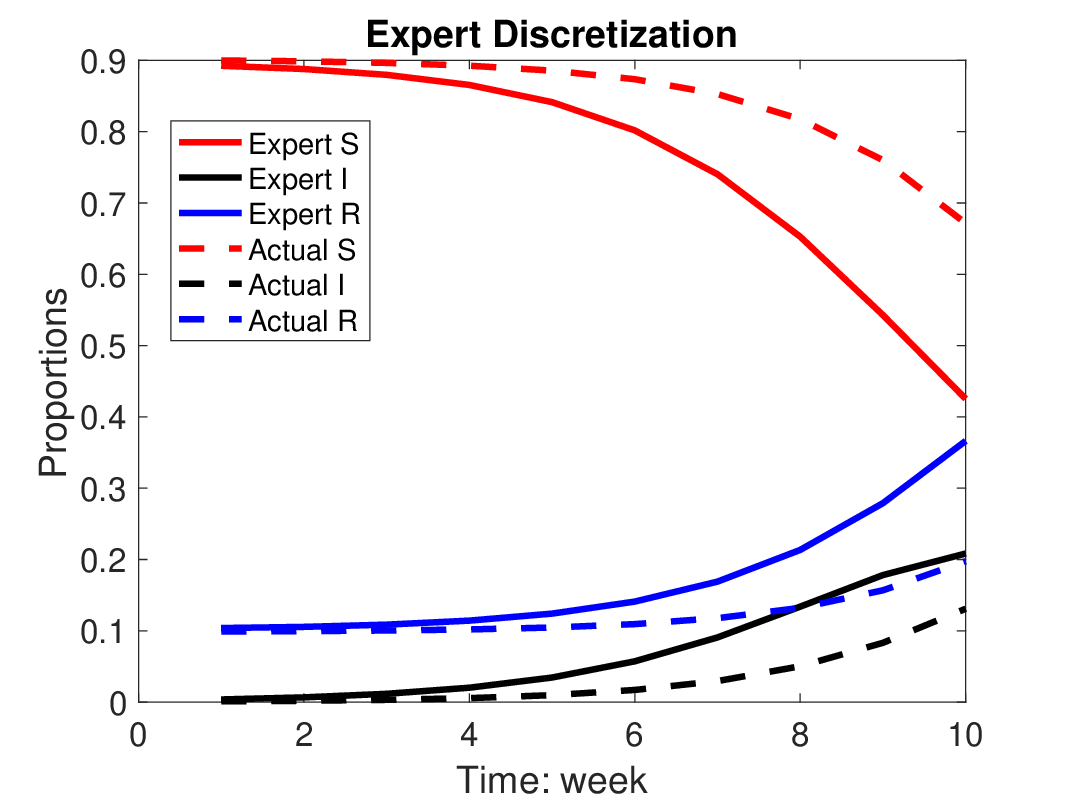}
 \caption{}
 \end{subfigure} %
 \begin{subfigure}{.49\textwidth}
 \centering
 \includegraphics[width=2.5in]{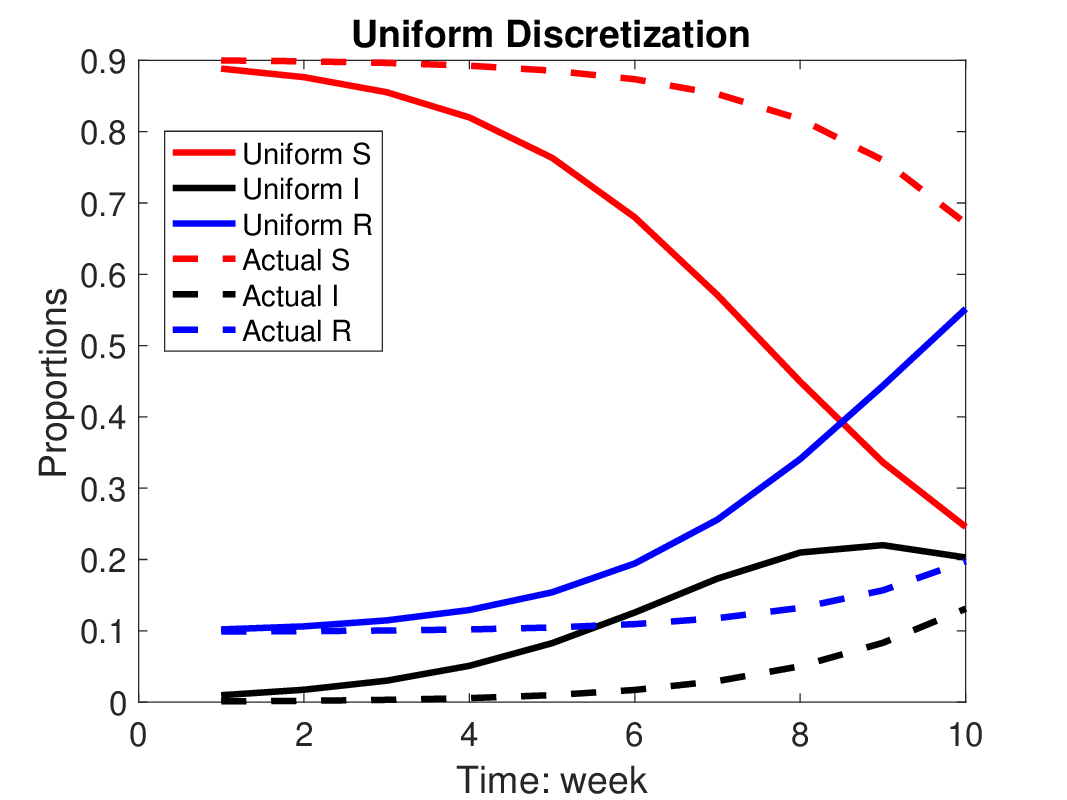}
 \caption{}
 \end{subfigure} %
 \caption{Comparison between trajectories generated from the GreedyCut discretization method against the uniform discretization method (using 300 discretization regions in total) given trajectories generated from SIR model. For each compartment S, I, and R, the GreedyCut discretization method can better capture the disease dynamics.  }
 \label{fig:figure3.2} 
\end{figure}

We use the same discretization levels (90, 150, 300, and 1200 discretization regions ) generated in the previous section for evaluation. To benchmark our model, we also generated uniform and InverseProportional discretization regions with the same number of discretization regions. Then, for all discretization methods, transition probability matrices were generated using Algorithm 2. 

We find that both the GreedyCut and InverseProportional discretization method is able to better approximate the disease dynamics over the uniform discretization method. As shown in Figure~\ref{fig:figure3.2}, a comparison among the GreedyCut, the InverseProportional, the expert, and the uniform discretization methods based on 300 discretization regions shows that the Markovian trajectories for both the GreedyCut and the InverseProportional discretization methods are closely aligned with the actual trajectories. However, the uniform discretization method shows poor approximation, especially showing incorrect trends for the proportion of the population infected over time -- where the proportion of the population infected over time starts to decline after week 8 in the Markovian trajectory, whereas the proportion of infected people over time increases in the entire time horizon in the actual trajectory. Additionally, the uniform overestimates the proportion of recovered populations by more than twice compared with the actual proportion of infected people.

\begin{table}[ht]
 \footnotesize
 \centering
 \begin{tabular}{|c|c|c|c|c|}
 \hline
 $B$ & GreedyCut [95\% CI]& InverseProportional [95\% CI]& Expert [95\% CI]&Uniform [95\% CI] \\
 \hline
 90&$0.1261 \quad[0.1100,0.1423]$&$0.1344 \quad[0.1160,0.1528]$& $0.1437\quad [0.1173,0.1700]$&$0.2399\quad[0.1627,0.3170]$\\
 \hline
 150&$0.1165\quad[0.1013,0.1318]$&$0.1117 \quad[0.0960,0.1274]$& $0.1301\quad [0.1143,0.1459]$&$0.1716\quad[0.1218,0.2215]$\\
 \hline
 300&$0.1088\quad[0.0943,0.1232]$&$0.1070 \quad[0.0927,0.1214]$& $0.1214\quad[0.1083,0.1346]$&$0.1326\quad[0.1118,0.1535]$\\
 \hline
 1200&$0.1071\quad[0.0935,0.1225]$&$0.1076 \quad[0.0929,0.1223]$& $0.1178\quad [0.1044,0.1312]$&$0.1197\quad[0.1059,0.1136]$\\ 
 \hline 
 \end{tabular}
 \caption{Mean squared error for the trajectories given different numbers of discretization regions}
 \label{tab:transitions_mse_uniform_greedy}
\end{table}

For all comparison pairs, the GreedyCut discretization method outperforms both expert and uniform discretization methods in the squared error between Markovian and actual trajectories. GreedyCut and InverseProportional discretization methods demonstrate similar performance overall comparison pairs. In Table~\ref{tab:transitions_mse_uniform_greedy}, we show the result of the comparison of $\mathbb{E}_{\hat{\theta}\in \hat{\Theta}}[||\{\tilde{\mathbf{X}}_t\}-\{\mathbf{X}_t\}||^2]$ over 100 samples and 10 time periods between the GreedyCut and uniform discretization methods. Both algorithms are able to improve the result of approximation when the number of discretization regions used increases, as expected. However, the improvement in approximations is small when the number of discretization regions is large, which suggests a high budget may not be necessary. When the number of allowable discretization regions is small due to the computational budget, the GreedyCut discretization method can provide a much better approximation than the uniform discretization method, and adding discretization regions may not add much accuracy.

Although GreedyCut and InverseProportional demonstrate similar MSE in estimating the actual trajectories, GreedyCut is still able to generate better MDP solutions. One possible reason could be that GreedyCut is finding a more accurate and relevant representation of the state space by minimizing $\mathbf{E}[||\{\bar{\mathbf{X}}_t\}-\mathbf{X}_t||^2]$, leading to better policy generation in the MDP. In contrast, the InverseProportional method primarily focuses on giving higher resolution to more frequently visited states, which results in good trajectory estimations but not necessarily in optimal MDP solutions. Additionally, the quality of InverseProportional discretization method is also sensitive to time horizon. For example, when the time horizon is long and the proportion of the population infected remains low for most of the time, the InverseProportional method will have an extremely high resolution in the region where the proportion of the population infected is low. However, this method will not provide sufficient resolution for other regions, potentially leading to suboptimal outcomes. This is also captured in the following example.

\subsection{Example 2: COVID-19}
COVID-19 led to a significant surge in infections within Los Angeles County (LAC). To mitigate the pandemic during its initial phases, LAC implemented a lockdown from the second week to the tenth week following March 1st, 2020, which marked the onset of the epidemic. In this example, we use an MDP with discretization regions to identify the optimal timing of imposing lockdowns in LAC to minimize the proportion of infected cases while considering the cost of a lockdown.

\subsubsection{Model Structure and Inputs}

To describe the disease dynamics of COVID-19 in LAC, we calibrated a SIR model that is stratified by health districts (HD) \citep{redelings2010years}, meaning that the model allows for heterogeneity in health outcomes across HDs. The transmission rates between HDs is also allowed to vary. The disease dynamics for HD $i$ are then described as follows:
 \begin{align*}      S^i_{t+1} & =S^i_{t} -\sum_{j}\beta_{ji}  S^i_t I^j_t \\
     I^i_{t+1}&= I^i_{t}+\sum_{j}\beta_{ji}  S^i_t I^j_t -\gamma I^i_t \\
     R^i_{t+1}&=R^i_{t}+\gamma I^i_t 
 \end{align*}
We use $\beta_{ij}$ to represent the transmission rate from HD $i$ to HD $j$, and all HDs are assumed to have the same clearance rate. We consider whether to implement a lockdown policy at each decision epoch, $\Delta t$, which has a duration of one week. Decisions need to be made over a total time horizon of 60 weeks (N=60). If lockdown is implemented, transmission will be reduced ($\beta$ decreases 80\%).

We assume there were 1000 infections (0.01\% of the total population) at the initial time epoch. This is consistent with the early stage of the COVID-19 epidemic in LAC where the proportion of the population infected remains a small proportion of the overall population. To calibrate the parameters of the stratified SIR model, we used empirical COVID-19 data of case counts to calibrate transmission rate $\beta$ and recovery rate $\gamma$ \citep{covidla}. LAC mobility data is also used to help capture the heterogeneity in transmission rate among HDs \citep{pems,yu2024extending}. 

We let the stage costs be the proportion of the population infected plus the dis-utility if lockdown is implemented. We assume that the dis-utility without a policy intervention is zero. However, determining the dis-utility associated with a lockdown is challenging. If the dis-utility is excessively low, the optimal choice consistently leans towards implementing the lockdown, which ignores the potential economic and social burden brought by the lockdown. On the contrary, if the dis-utility proves to be excessively burdensome, it will never be enforced. To better reflect this tradeoff, we assume the dis-utility of implementing a one-week lockdown is 0.005, implying it equates to the dis-utility of 0.5\% of the population infected per epoch.

We create three discretized MDPs with 150 discretization regions using the GreedyCut discretization method and compare outcomes against those of the InverseProportional, expert, and uniform discretization methods. This will guarantee the completion of model construction within an hour for all discretization methods. In this COVID-19 example, the brute force results cannot be generated within a reasonable time as there are $2^{60}$ different cases compared to only $2^{10}$ different cases in the previous example. Therefore, we directly compare the MDP objective values generated by using the identified policies on the true disease model to evaluate each algorithm's performance.

\subsubsection{MDP Results}\label{covidresults}

\begin{figure}[ht]
 \centering
 \begin{subfigure}{.49\textwidth}
 \centering \includegraphics[width=2.5in]{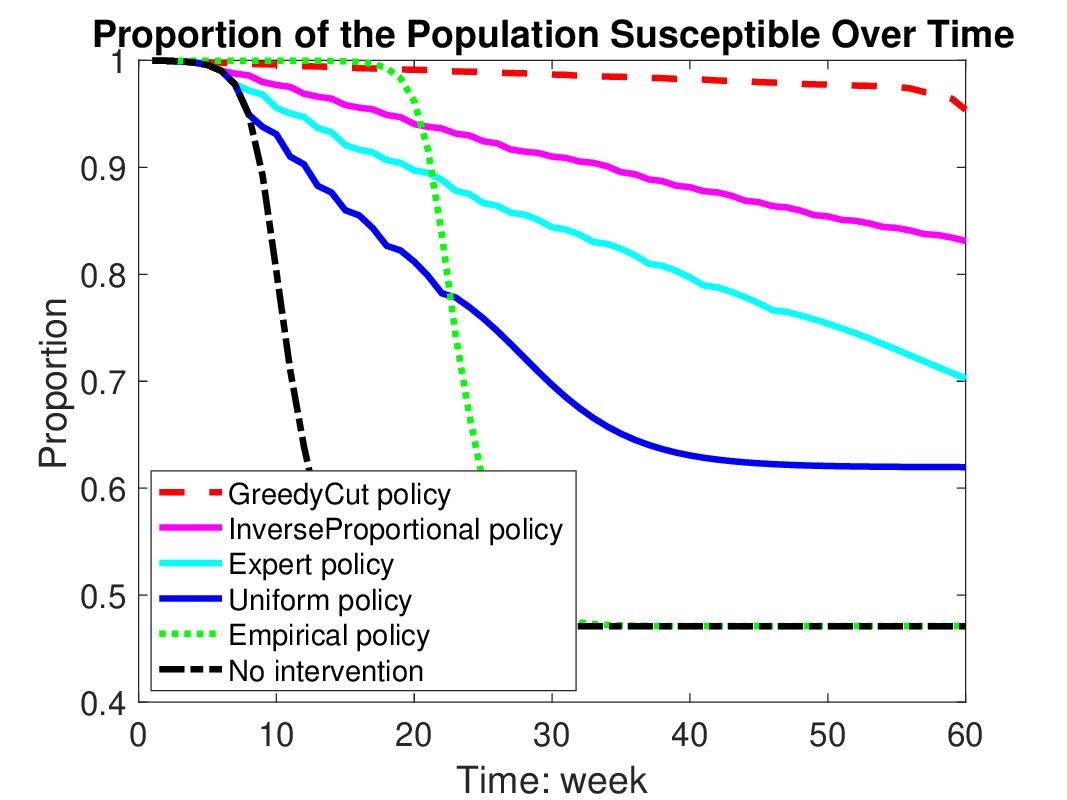}
 \caption{}\label{fig:fig_a}
 \end{subfigure} %
 \begin{subfigure}{.49\textwidth}
 \centering
 \includegraphics[width=2.5in]{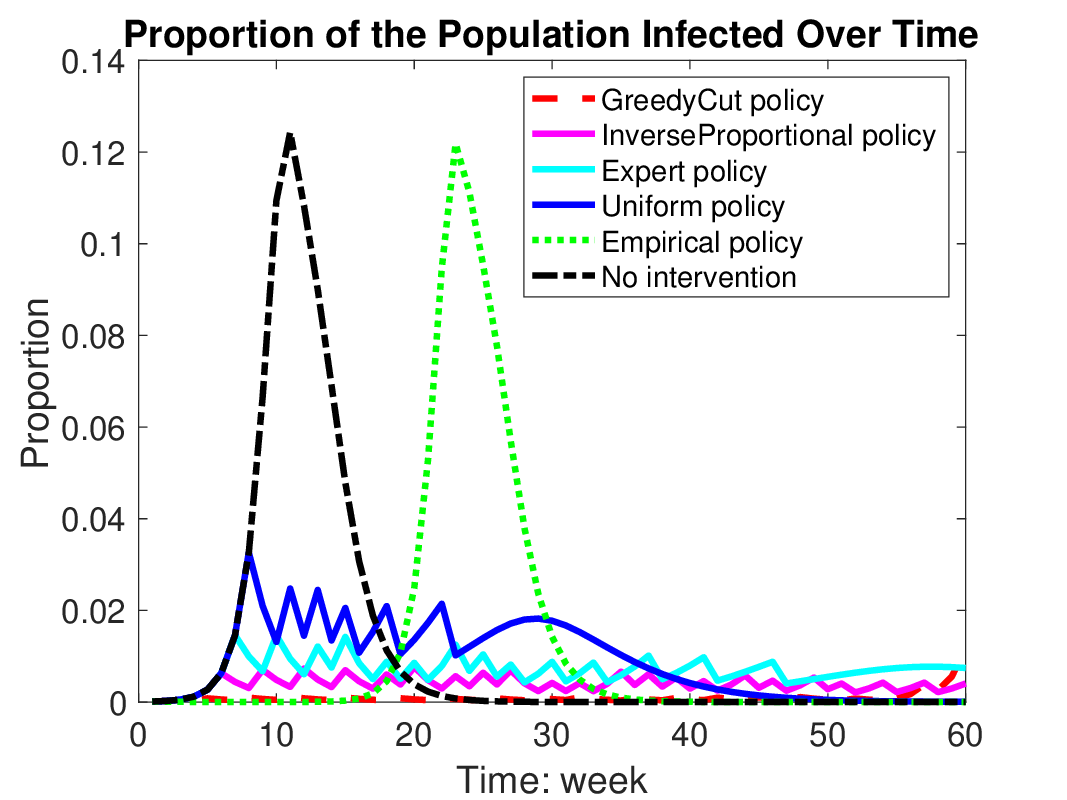}
 \caption{}
 \end{subfigure} %
 \caption{ Proportions of the population that is susceptible/infected over time.  }
 \label{fig:covid1} 
\end{figure}

We compared the optimal action recommended by the discretized state MDP from the GreedyCut, InverseProportional, expert, and uniform discretization methods. We find that our GreedyCut algorithm outperforms the other discretization methods by identifying a better MDP optimal solution with a smaller cumulative proportion of the population infected. Figure~\ref{fig:covid1} shows a comparison of disease dynamics across different policies. Compared with no lockdown, the empirical policy in LAC (lockdown from week 2 to week 10) does not prevent but rather postpones infections (the total cumulative proportion of the population infected over time is 0.7504 in the empirical policy, and 0.7508 if no intervention is used). All uniform, InverseProportional, expert, and GreedyCut discretized MDPs are able to reduce the cumulative proportion of infections and the peak of infections. The GreedyCut discretization method outperforms the uniform discretization method in terms of the overall reduction in the proportion of the population infected by 0.4793 over the 60-week time horizon.

\begin{figure}[ht]
 \centering 
 \includegraphics[width=4in]{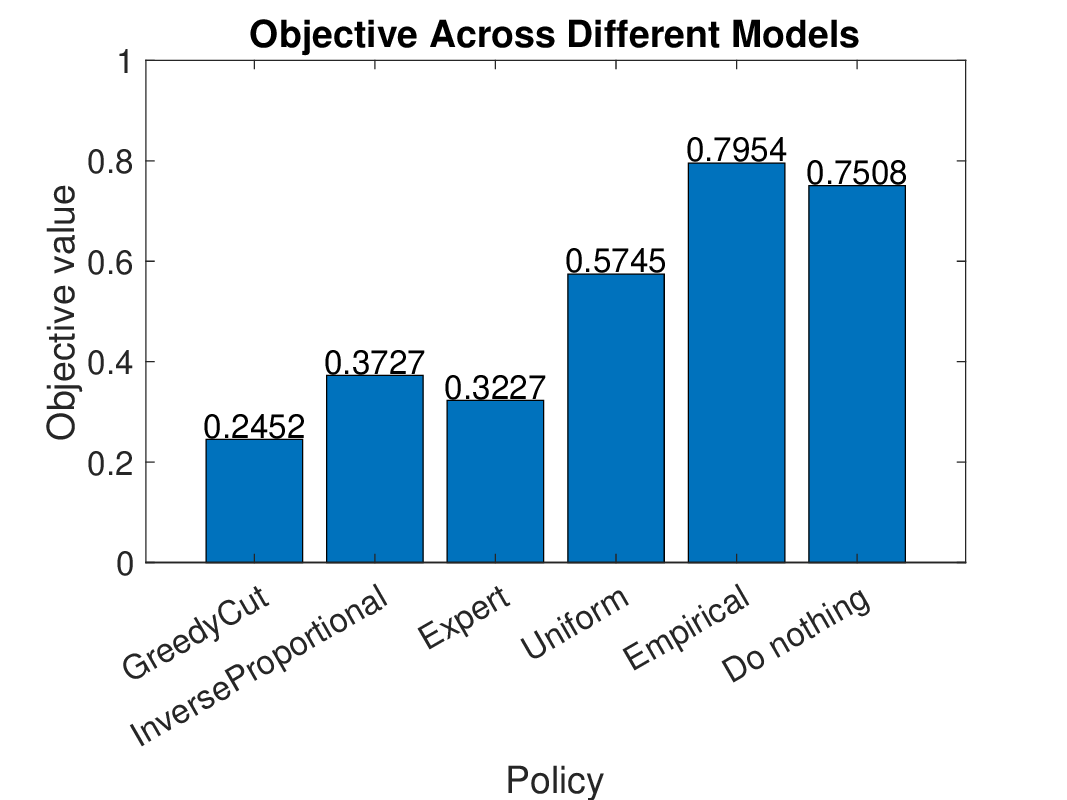} 
 \caption[Comparison between objectives.]{Comparison between objective values across policies from GreedyCut MDP, InverseProportional MDP, uniform MDP, empirical policy, and no policy.}
 \label{fig:obj} 
\end{figure}

GreedyCut outperforms uniform, InverseProportional, expert, empirical, and `do nothing' policies by achieving the lowest objective value (proportion of infected people over time plus lockdown disutility). Figure~\ref{fig:obj} compares the objective value among different models evaluated on the ground-truth disease (compartmental) model. The 8-week empirical lockdown policy LAC imposed has a lower objective value than doing nothing after considering the cost of lockdown, as it was not able to reduce infections while incurring lockdown disutility costs. Both uniform and GreedyCut MDPs provide a better solution. GreedyCut MDP is able to improve the objective value from doing nothing by 67\%, from InverseProportional by 34\%, and outperforms the uniform MDP outcome by 57\%.

This example demonstrates that the GreedyCut algorithm is able to provide a solution that has a smaller total cost compared to the empirical lockdown policy in LAC and the policy generated by the other discretized MDPs. Even when the number of discretization regions is limited for each compartment (for example, a stratified compartmental model includes death, hospitalizations, exposed, etc.), the GreedyCut discretization method can generate quality solutions with low total discounted costs. Moreover, we found that the uniform discretization method is not able to generate a near-optimal solution as its optimal value exceeds twice the objective value of the GreedyCut discretized MDP. 

\subsubsection{Sensitivity Analysis: Uncertainty in Transmission Rate}
We calibrated the transmission rate of a COVID-19 compartmental model using empirical COVID-19 data and traffic data as a proxy for the contact matrix. However, these rates may be subject to measurement uncertainty. We therefore perform sensitivity analysis where the discretization is based on the original calibrated values, but the actual transmission rate is (1) 30\% lower and (2) 30\% higher to understand the performance of the discretization methods when input parameters are erroneous. The MDP solutions are generated using the actual transmission rate, while the discretization uses the original calibrated value. This allows us to compare the robustness and reliability of different discretization regions. In each case, we found that GreedyCut was still able to generate the best policy with the lowest objective value (details shown in Appendix Figure~\ref{fig:senstivity}). This demonstrates that discretization regions generated by GreedyCut are robust to variations in the disease variables.

\subsubsection{Extension: MDP with at Most Two Policy Switches }
\begin{figure}[ht]
 \centering 
 \includegraphics[width=6in]{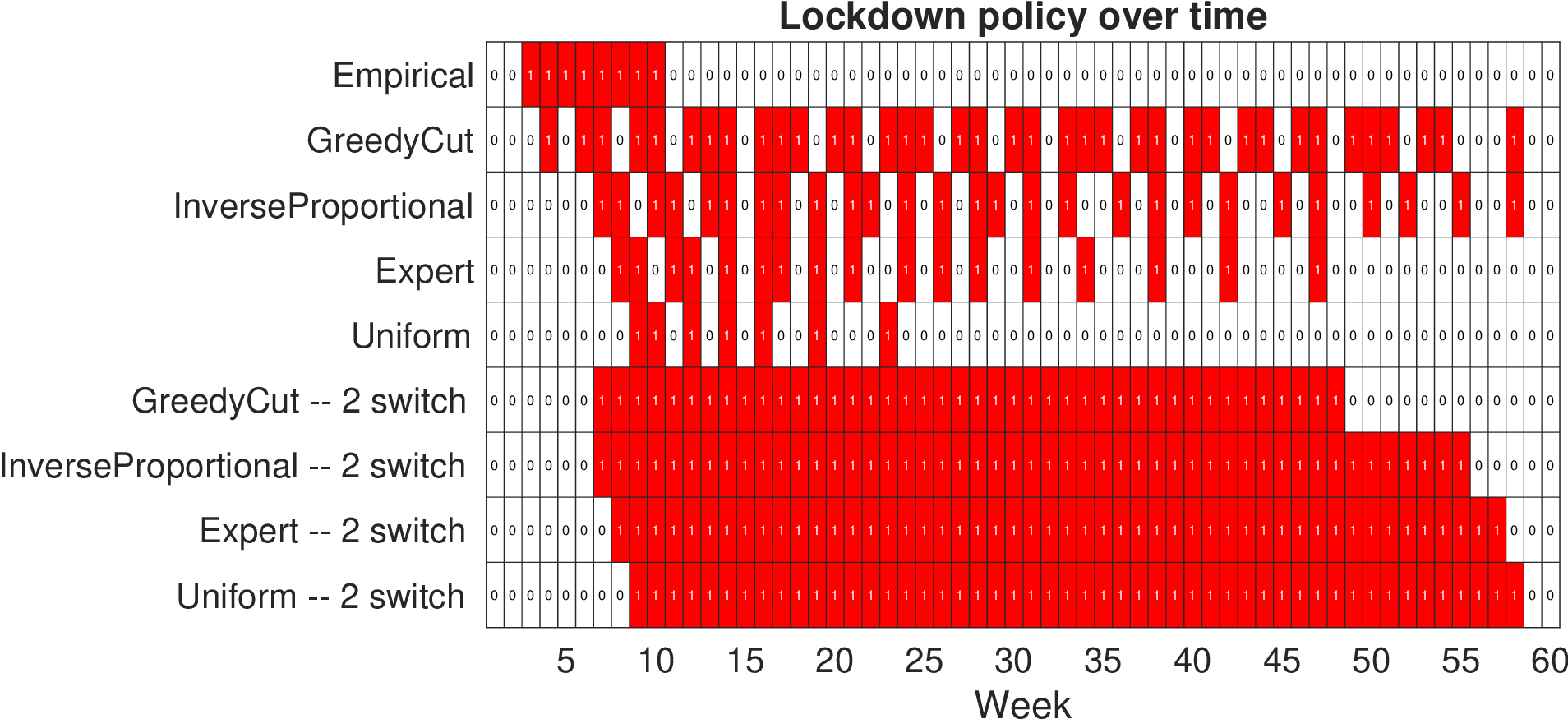} 
 \caption{Lockdown policy. The GreedyCut discretized MDP recommends starting the lockdown on week 7 for 42 weeks. The InverseProportional discretized MDP recommends starting the lockdown on week 7 for 49 weeks. The expert discretized MDP recommends starting the lockdown on week 8 for 50 weeks. The uniform discretized MDP recommends starting the lockdown on week 9 for 50 weeks.  }
 \label{fig:policy_seq} 
\end{figure}

In Section~\ref{covidresults}, all discretized MDPs recommended a policy with many policy switches where lockdown would be imposed for many short durations. For example, the policy from the uniform discretized MDP recommends lockdown every few weeks in weeks 9-23 (shown in Figure~\ref{fig:policy_seq}). This is not practical, as inconsistency in policies can lead to poor adherence or even psychological issues \citep{webster2020improve,pedrozo2020perceived}. In this section, we consider the same problem with additional constraints where we allow the policy to switch at most twice (once from no lockdown to lockdown, and once from lockdown to no lockdown). We set up the COVID-19 dynamics in the same way as in  Section~\ref{covidresults}. 

\begin{figure}[ht]
 \centering
 \begin{subfigure}{.49\textwidth}
 \centering \includegraphics[width=2.5in]{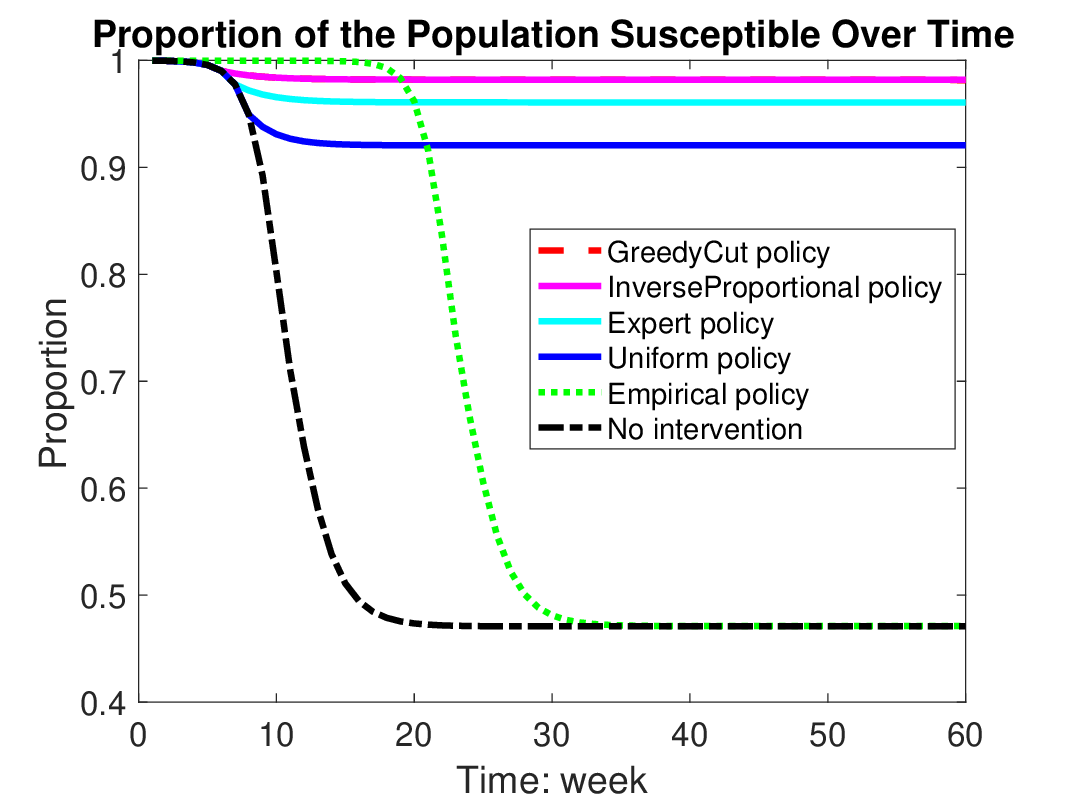}
 \caption{}\label{fig:fig_a}
 \end{subfigure} %
 \begin{subfigure}{.49\textwidth}
 \centering
\includegraphics[width=2.5in]{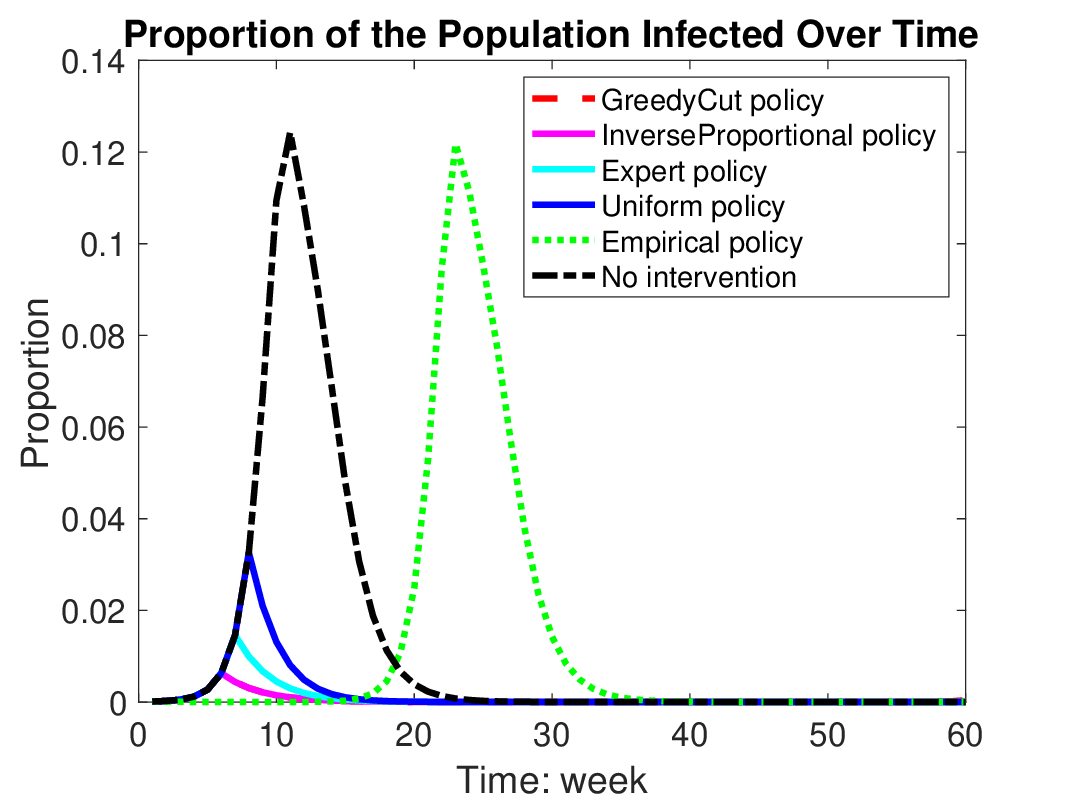}
 \caption{}
 \end{subfigure} %
 \caption{ Proportion of the population susceptible/infected over time for different MDPs with policy constraints.}
 \label{fig:covid1_switch} 
\end{figure}

With at most two policy switches (one lockdown duration), the GreedyCut discretized MDP recommends a shorter lockdown duration than the uniform discretized MDP and an earlier lockdown initiation date. Figure~\ref{fig:policy_seq} shows the lockdown policy outcomes given the disease dynamics of the ground-truth model. The GreedyCut discretized MDP recommends starting the lockdown on week 7 for a duration of 42 weeks. The InverseProportional discretized MDP recommends starting the lockdown at the same time as GreedyCut but with a longer lockdown duration. The expert and uniform discretized MDP recommend initiating the lockdown in weeks 8 and 9, lasting for a duration of 50 weeks. Due to the highly transmissible nature of COVID-19, a lockdown of over 40 weeks is needed to reduce transmission.  Figure~\ref{fig:covid1_switch} compares the trajectories using policies from the uniform discretized MDP and GreedyCut discretized MDP. The GreedyCut discretized MDP is able to generate a policy that reduces the cumulative proportion of the population infected by 0.0864 while using fewer weeks of lockdowns compared to the uniform discretized MDP. Additionally, we observed that most countries do not impose prolonged lockdowns. Therefore, we considered an 11-week lockdown using the different discretized MDPs. The GreedyCut discretized MDP recommends starting the lockdown from week 12, resulting in a cumulative infection proportion of 0.458. The uniform, expert, and InverseProportional discretized MDPs recommend starting the lockdown in weeks 10, 13, and 16, respectively, resulting in cumulative infection proportions of 0.5834, 0.5399, and 0.7012. All these values are higher than those achieved by implementing longer lockdown policies, confirming the necessity of imposing longer lockdowns.

With the additional constraint on the policy switches, the GreedyCut discretized MDP consistently generates a better solution than other discretized MDPs and other policies we considered in this analysis. In addition, the GreedyCut discretized MDP recommends a shorter lockdown duration compared with other discretized MDPs which could reduce the economic and societal burden brought by the lockdown.

\section{Conclusions}\label{conclusions_discretization}

In this paper, we introduce a novel algorithm for formulating an MDP framework tailored for continuous or large state-space problems in repeated decision-making for infected disease control. In our numerical analyses, we found that our algorithm provides better MDP solutions than the other discretized MDPs for the models we evaluated. Our approach better approximates the true value function than the uniform discretized MDP, therefore leading to a better policy with a lower optimality gap. Compared to other discretization regions, our method demonstrates better performance across different discretization budgets, particularly showing notable benefits when the number of discretization regions is small. This may be particularly pertinent for a compartmental model with many compartments, as the resultant number of discretization regions for each compartment may be extremely limited. In this case, a uniform discretization method may result in a poor estimation of disease dynamics.

We found that our algorithm is able to provide a better state discretization than the uniform discretization method in approximating disease dynamics for the examples considered. Our approach generates smaller regions for states with a higher likelihood of being visited, while also preserving some regions for states with a lower but still possible likelihood, leading to a robust and reliable discretization. Using the GreedyCut method would substantially improve the approximation quality, thus resulting in an improved decision-making process.

%Additionally, we offer an effective method to compute the transition probability matrices for formulating MDPs given $f(\mathbf{X}_t,\pi_t)$, which may be a compartmental or simulation-based model. This helps us incorporate other common disease modeling frameworks into the MDP decision-making framework, facilitating the seamless integration of diverse healthcare applications. This approach is not only straightforward to implement but also effectively captures the complex function $f(\mathbf{X}_t,\pi_t)$ that represents the dynamics of the disease in transition probability matrices.

With a small number of discretization regions, the time spent discretizing the state using our approach is considerably smaller than the time required to produce transition probability matrices. In our examples, the time needed to formulate the discretized MDP using our approach is nearly equivalent to the time necessary for the construction of a uniform discretized MDP or other benchmark methods. Therefore, our approach could offer an improvement to the MDP solution without substantially increasing computational expense under limited budgets.

We provide numeric examples to demonstrate the efficiency and effectiveness of our algorithm in discretizing continuous-state decision-making problems. Benchmarking the performance with the different discretization methods, we demonstrate that our algorithm is able to generate preferable discretization regions with a limited budget that is a good proxy of the ground truth system. We also demonstrate that our algorithm can generate better policies in both synthetic examples and a COVID-19 example. In the synthetic example, our algorithm outperforms other discretization regions in all metrics for different discretization budgets. In the COVID-19 example, our algorithm improves the objective by nearly 100\% from the uniform discretization.

Our numerical analysis also generated policy implications for social distancing policy during COVID-19 in LAC. The first policy implication is that the threshold of implementing the lockdown depends on the proportions of susceptible and infected (recommending implementing the lockdown if the proportions of infected and susceptible are above certain numbers). When the proportion of the population infected increases, the threshold of implementing the lockdown on the proportion of the susceptible population decreases. This is because a less susceptible population is needed to spread the disease as the infected population grows. Similarly, when the proportion of the susceptible population increases, the threshold of implementing the lockdown on the proportion of the infected population decreases. Secondly, a short lockdown interval may not effectively reduce the total number of cases, instead only delays the epidemic peak. To more effectively control cases, a prolonged lockdown period is needed.

We must acknowledge several limitations of this work. The GreedyCut algorithm may not find the discretization that globally minimizes the cost function. The performance gap between the GreedyCut and the uniform discretization methods is small with a large discretization budget. The GreedyCut algorithm may have computational difficulty if there is a large action space and does not consider continuous action spaces; this leaves an interesting optimization direction for future studies. The output of the GreedyCut algorithm may be sensitive to the choice of cost function; different choices may result in widely different discretization choices, thus requiring reevaluation of the objective function.

Despite these limitations, we believe that this work provides an effective and easy-to-handle scheme for dealing with decision-making problems in large or continuous state spaces. Our paper provides insight into future work on improving the discretization of solving large-scale MDPs.

\section*{Acknowledgements}
Financial support for this study was provided in part by a National Science Foundation (grant no. 2237959).

\bibliographystyle{apacite}
\bibliography{interactapasample}

\newpage
\setcounter{section}{1}
\section*{Appendix}\label{appendix}
\renewcommand{\thetable}{Appendix \arabic{table}}
\renewcommand{\thefigure}{Appendix \arabic{figure}}
\setcounter{table}{0} % This sets the table number to 2 (since it increments afterward)
\setcounter{figure}{0} % This sets the table number to 2 (since it increments afterward)
\subsection{SAA Algorithm for Constructing Transition Matrices}
\begin{algorithm}
\caption{Generating Transition matrix}\label{sample}
\begin{algorithmic}[1]
\Procedure{Generate}{$f(\mathbf{X}_t,\pi_t),G$}\Comment{$f(\mathbf{X}_t,\pi_t)$ is the ground-truth discrete time model, $G$ is the discretization }
\State Let $P$ be a $|\bar{\mathcal{X}}|\times |\bar{\mathcal{X}}|\times |A|$ transition matrix with all zeros and each state represents a discretized state from the discretized state space $\bar{\mathcal{X}}$ defined by $G$
\For {each policy intervention $\pi_t\in A$}
\For {each discretized state $i$ from $\bar{\mathcal{X}}$}
\State Uniformly draw $c$ number of samples ($\hat{\mathbf{X}}_0$) within the region that contains $i$ (including a centroid in this region)
\For {each sample $\hat{\mathbf{X}}_0$}
\State Compute $\hat{\mathbf{X}}_1 = f(\hat{\mathbf{X}}_0,\pi_t)$
\State Find the discretized state $j$ such that the discretized region containing $j$ also contains $\hat{\mathbf{X}}_1$
\State $P(j|i,\pi_t) = P(j|i,\pi_t)+1$
\EndFor
\EndFor
\State Normalize $P$ to make it a stochastic matrix
\EndFor
\State \Return{$P$}
\EndProcedure
\end{algorithmic}
\end{algorithm}
\subsection{Total Runtime by Different Algorithms}
\begin{table}[h]
    \centering
    \begin{tabular}{c|c|c|c|c|}
         &Uniform&GreedyCut&InverseProportional&ExpertOpinion  \\
         \hline
         90&434&436&436&434\\ 
         \hline
         150&1226&1236&1228&1226\\ 
                  \hline

         300&3625&3660&3630&3625\\ 
                  \hline

         1200&10,560&10,721&10,568&10,560\\
                  \hline

    \end{tabular}
    \caption{Total runtime by different algorithms, including discretization, generation of the transition probability matrix, and solving the MDP (in seconds).}
    \label{tab:total_runtime}
\end{table}
\subsection{MDP Solutions of GreedyCut Using Different Cost Functions}
\begin{table}[ht]
 \scriptsize
 \centering
 \begin{tabular}{|c|c|c|c|c|c|c|c|c|}
 \hline
  $B$ &  \multicolumn{2}{c|}{ACC} & \multicolumn{2}{c|}{MSE}\\
  \hline
   & GreedyCut (Sum of Squared Error) & GreedyCut (MAE)&GreedyCut (Sum of Squared Error) & GreedyCut (MAE) \\
   \hline
   90&0.9657&0.9457&4.3239e-04&8.5366e-04\\
   \hline
   150&0.9850&0.9357&1.7896e-04 & 6.6263e-04\\
   \hline
   300&0.9787&0.9733&6.5032e-05&1.7771e-04\\
   \hline
   1200&0.9893&0.9900&3.0731e-06 &1.8360e-05\\
   \hline
  $B$ &  \multicolumn{2}{c|}{E2} & \multicolumn{2}{c|}{Opt. Gap}\\
  \hline
   & GreedyCut (Sum of Squared Error) & GreedyCut (MAE)&GreedyCut (Sum of Squared Error) & GreedyCut (MAE) \\
   \hline
   90&0.0689 &0.0905&0.0033&0.0060\\
   \hline
   150&0.0435& 0.0839&0.0011&0.0089\\
   \hline
   300&0.0233&0.0457&0.0018 &0.0017\\
   \hline
   1200&0.0048&0.0126&4.5576e-04 &4.1921e-04\\
   \hline
 \end{tabular}
 \caption{Comparison on MDP solutions between sum of the squared error and mean absolute error}
 \label{tab:appendix_tab1}
\end{table}
\newpage
\subsection{Matrix Convergence Over Number of Samples}
\begin{figure}[h]
    \centering
    \includegraphics[width=0.5\linewidth]{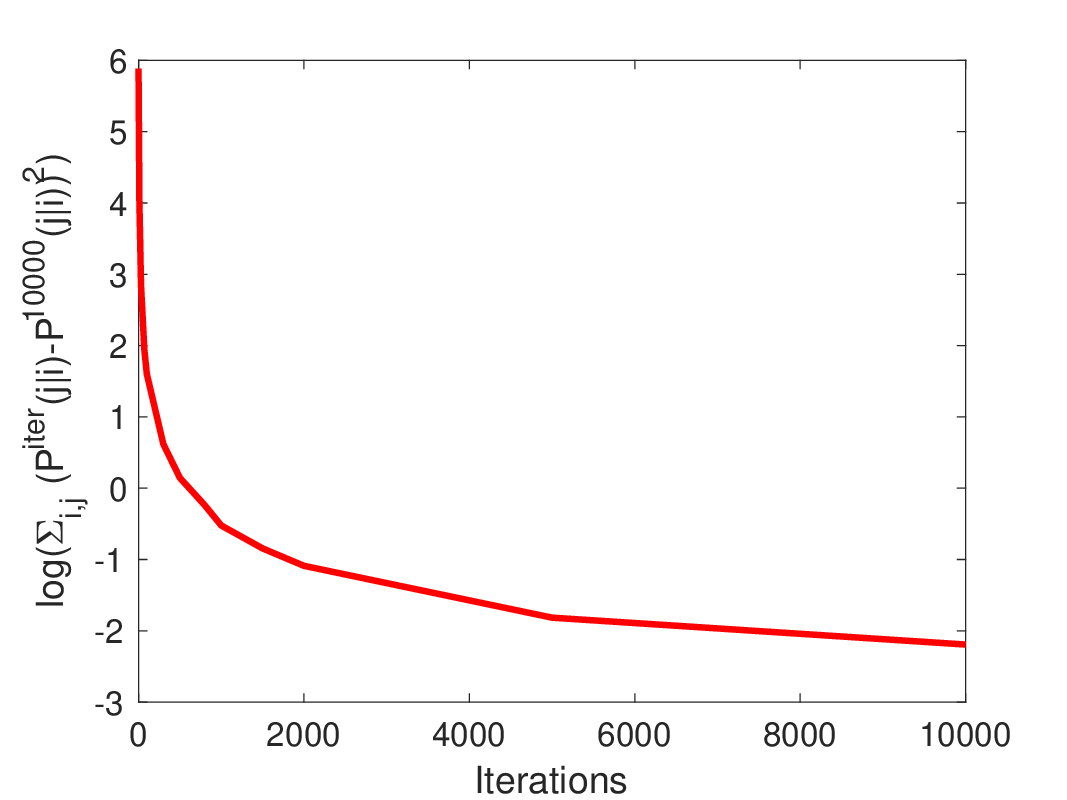}
    \caption{Matrix Convergence Over Time}
    \label{fig:convergence}
\end{figure}
\newpage
\subsection{Distribution of Discretization Regions Generated}

\begin{figure}[ht]
    \centering
    \includegraphics[width=0.5\linewidth]{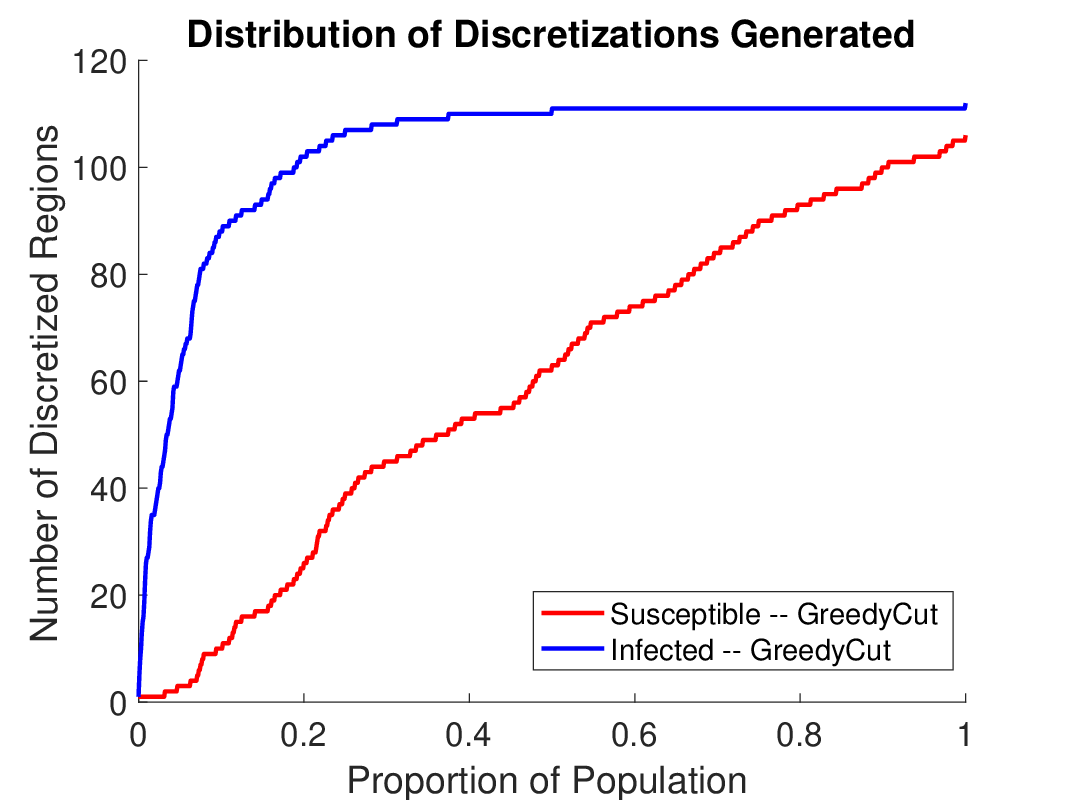}
    \caption{The GreedyCut algorithm allocated most of its cuts to states that are more likely to be visited based on sample trajectories. Unlike the InverseProportional method, which focuses solely on frequently visited states, GreedyCut also made cuts in less frequently visited but still possible regions, potentially offering a more reliable and robust discretization.}
    \label{fig:distribution_cuts}
\end{figure}
\subsection{MDP Outcomes Using Different Numbers of State Time Pairs}
\begin{table}[h!]
 \scriptsize
 \centering
 \begin{tabular}{|c|c|c|c|c|}
 \hline
   $B$&\multicolumn{2}{c|}{ACC}& \multicolumn{2}{c|}{MSE}\\
\hline
& GreedyCut (3000 pairs)&GreedyCut (6000 pairs) &GreedyCut (3000 pairs)&GreedyCut (6000 pairs) \\
   \hline
   90&0.9657&0.9597&4.3239e-04&3.1057e-04\\
   \hline
   150&0.9850&0.9711&1.7896e-04&1.3284e-04\\
   \hline
   300&0.9787& 0.9772&6.5032e-05& 5.4063e-05\\
   \hline
   1200&0.9893&0.9927&3.0731e-06&2.3932e-06\\
   \hline
  \hline
  $B$&\multicolumn{2}{c|}{E2}& \multicolumn{2}{c|}{Opt. Gap}\\
\hline
& GreedyCut (3000 pairs)&GreedyCut (6000 pairs) &GreedyCut (3000 pairs)&GreedyCut (6000 pairs) \\
 \hline
90&0.0689&0.0546&0.0033& 0.0036\\
 \hline
150&0.0435&0.0350&0.0011&0.0019\\
 \hline
300&0.0233& 0.0212&0.0018& 0.0019\\
 \hline
 1200&0.0048&0.0041&4.5576e-04& 2.8458e-04\\ 
 \hline 
 \end{tabular}
 \caption{Comparison on MDP solutions}
 \label{tab:mdp_comapre_state_pairs}
\end{table}
\newpage
\subsection{Sensitivity Analysis: Uncertainty in Transmission Rate}
\begin{figure}[ht]
 \centering

 \begin{subfigure}{.49\textwidth}
 \centering
 \includegraphics[width=3in]{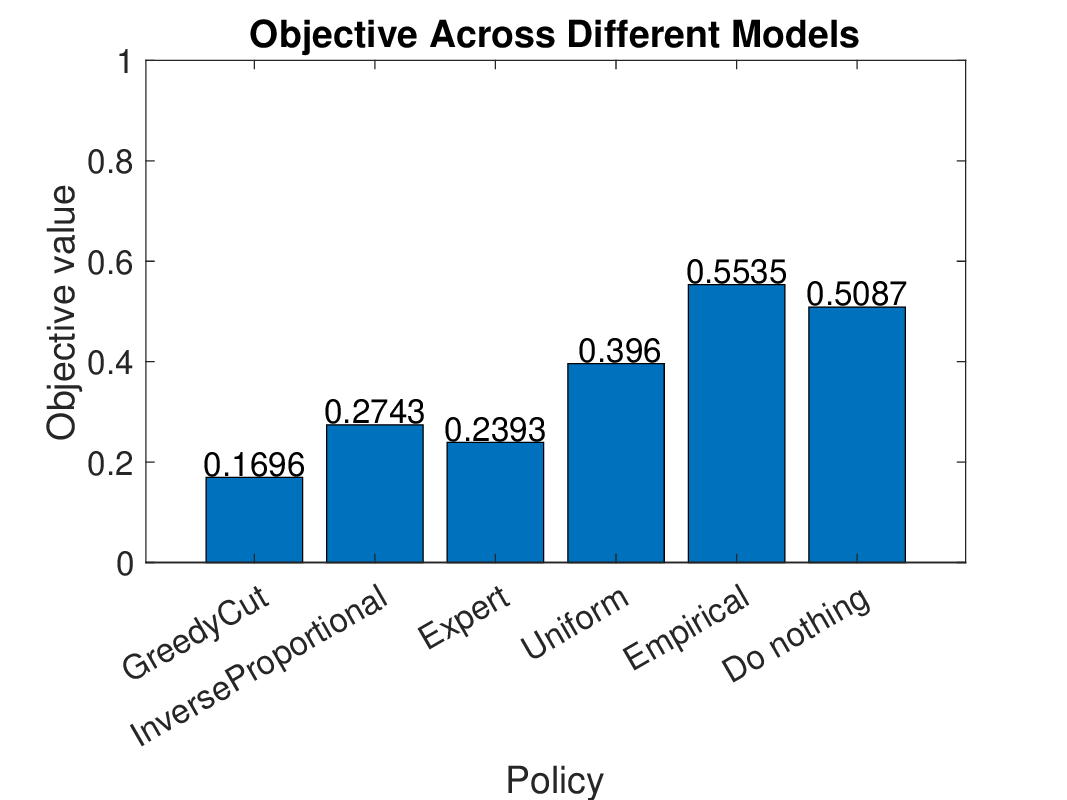}
 \caption{}
 \label{fig:30lower} 
 \end{subfigure} %
 \begin{subfigure}{.49\textwidth}
 \centering
 \includegraphics[width=3in]{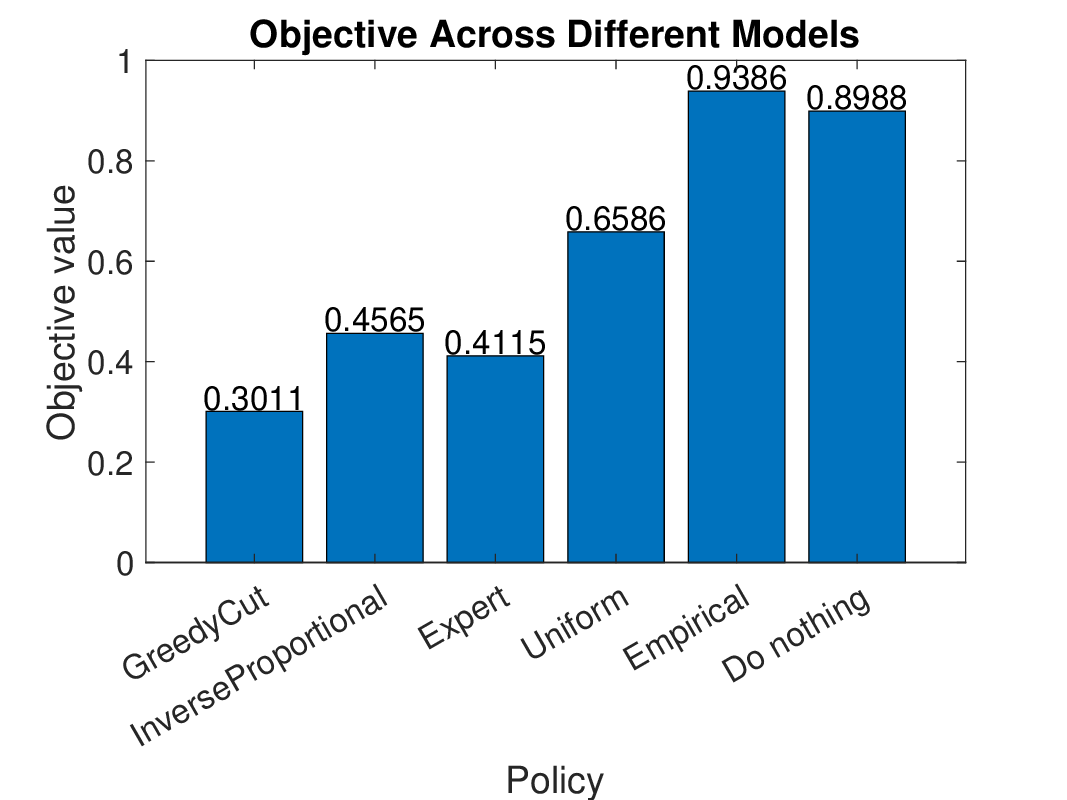}
 \caption{}
 \label{fig:30more} 
 \end{subfigure} %
 \caption{The objective of MDP solutions. (1) the actual transmission rate is 30\% less than the calibrated values used in the discretization generation. (2) the actual transmission rate is 30\% more than the calibrated values used in the discretization generation.}
 \label{fig:senstivity} 
\end{figure}

\end{document}